\documentclass[twoside]{article}
\usepackage{graphicx} % Required for inserting images

\usepackage[english]{babel}
\usepackage{csquotes}
\usepackage[toc,page]{appendix}

\usepackage[
  a4paper,
  inner=2.2cm, outer=2.6cm,
  top=3cm,    bottom=3cm,
  bindingoffset=0.5cm,
  headsep=1em
]{geometry}
 
\usepackage{amsmath,amsthm,amssymb}
\usepackage{mathtools}
\usetagform{default}

 \usepackage{mathrsfs}
 \usepackage{esint}
\usepackage[nopatch=eqnum]{microtype}

\usepackage{mathspec}
\setallmainfonts[
  ItalicFont     = OldStandard-Italic.otf,
  BoldFont       = OldStandard-Bold.otf,
  BoldItalicFont = OldStandard-BoldItalic.otf
]{OldStandard-Regular.otf}

\usepackage{graphicx}
\usepackage{xcolor}
\usepackage{booktabs}
\usepackage{array}
\usepackage{caption}
\usepackage{subcaption}
\usepackage{tikz}
\usepackage{enumitem}
\usepackage{bm}
\usepackage{xfrac}
\usetikzlibrary{calc,arrows.meta,decorations.pathmorphing}
\usepackage{csquotes}

\usepackage{titlesec}

\titleformat{\section}
  {\normalfont\large\scshape\filcenter}
  {\thesection.}{1em}{}

\titlespacing*{\section}
  {0pt}{2.5ex plus 1ex minus .2ex}{1.5ex plus .2ex}

\titleformat{\subsection}[runin]
  {\normalfont\bfseries}
  {\thesubsection.}{0.5em}{}[.]

\titlespacing*{\subsection}
  {0pt}{1.5ex plus 0.5ex minus 0.2ex}{1em}

\titleformat{\subsubsection}[runin]
  {\normalfont\itshape}
  {\thesubsubsection.}{0.5em}{}[.]

\titlespacing*{\subsubsection}
  {0pt}{1.2ex plus 0.4ex minus 0.2ex}{1em}

\titleformat{\chapter}[display]
  {\normalfont\huge\bfseries}
  {\chaptertitlename\ \thechapter}
  {20pt}
  {\Huge}

\titlespacing*{\chapter}{0pt}{50pt}{80pt}
 
\usepackage[
  colorlinks=true,
  linkcolor=black,
  citecolor=black,
  urlcolor=blue!60!black,
  bookmarksnumbered=true
]{hyperref}
\usepackage[nameinlink,capitalise,noabbrev]{cleveref}

\usepackage{fancyhdr}
\numberwithin{equation}{section}

\usepackage{aliascnt}

\newcommand{\newshared}[3]{%
  \newaliascnt{#1}{thm}%
  \newtheorem{#1}[#1]{#2}%   <- [#1] aggiunto
  \aliascntresetthe{#1}%
  \crefname{#1}{#2}{#3}%
  \Crefname{#1}{#2}{#3}%
}

\theoremstyle{plain}
\newtheorem{introthm}{Theorem}

\crefname{introthm}{Theorem}{Theorems}
\Crefname{introthm}{Theorem}{Theorems}

\newtheorem{thm}{Theorem}[section]
\crefname{thm}{Theorem}{Theorems}
\Crefname{thm}{Theorem}{Theorems}

\newshared{prop}{Proposition}{Propositions}
\newshared{lem}{Lemma}{Lemmas}
\newshared{cor}{Corollary}{Corollaries}

\newtheorem*{thm*}{Theorem}   % unnumbered, no aliascnt needed
\newtheorem*{cor*}{Corollary}

\theoremstyle{definition}
\newshared{dfn}{Definition}{Definitions}
\newshared{ex}{Example}{Examples}
\newshared{prb}{Problem}{Problems}
\newshared{ass}{Assumption}{Assumptions}

\theoremstyle{remark}
\newshared{rmk}{Remark}{Remarks}
\newshared{notation}{Notation}{Notations}
 
\renewcommand{\div}{\operatorname{div}}

\newcommand{\tr}{\operatorname{tr}}
\newcommand{\spt}{\operatorname{spt}}
\newcommand{\dist}{\operatorname{dist}}

\newcommand{\id}{\operatorname{id}}
\newcommand{\lip}{\operatorname{Lip}}
\newcommand{\loc}{\mathrm{loc}}

\newcommand{\R}{\mathbb{R}}

\renewcommand{\S}{\mathbb{S}}
 
\newcommand{\Mink}{\R^{m,1}}
\newcommand{\II}{\mathrm{II}}
\newcommand{\norm}[1]{\lVert#1\rVert}

\newcommand{\set}[1]{\{#1\}}

\newcommand{\X}{\mathfrak{X}}
\newcommand{\eps}{\varepsilon}

\newcommand{\Haus}{\mathscr{H}}

\newcommand{\babla}{\overline\nabla}

\newcommand{\de}{\mathrm d}
\newcommand{\del}{\partial}

\title{\textbf{Spacelike graphs with prescribed mean curvature and capillary boundary datum}}
\author{%
  \textsc{Lorenzo Maniscalco}
}
\date{}

\begin{document}

\maketitle

% \scriptsize \begin{center} Dipartimento di Matematica ``Giuseppe Peano'' \\ 
%  Università degli Studi di Torino, Via Carlo Alberto 10, 10123 Torino, Italy\\
% \texttt{lorenzo.maniscalco@unito.it}
% \end{center}
% \bigskip

\begin{abstract}
    \noindent We study space-like graphs in Minkowski space over bounded domains in $\mathbb{R}^m$, with prescribed mean curvature $\rho$ and capillary boundary datum $\psi$. The corresponding equation is the Euler--Lagrange equation of a convex, nonsmooth functional. On any bounded convex domain $\Omega$, we prove that the graph of the minimiser contains no interior light segments whenever $\rho\in L^{m-1}(\Omega)$ and $\psi\in L^{m-1}(\partial\Omega)$ satisfy the compatibility condition. On a ball $B$, we also consider measures $\rho$ whose singular part is supported on a compact set $\Sigma$ with $\mathcal{H}^1(\Sigma)=0$ and whose absolutely continuous part has a locally square-integrable density away from $\Sigma$. For compatible data with $\psi\in L^q(\partial B)$ for some $q>1$, the minimiser solves the capillary problem weakly up to the boundary away from the closure of the union of its light segments.

    % \bigskip

    \par
\end{abstract}

\bigskip
\medskip
% \noindent\footnotesize\textbf{2020 Mathematics Subject Classification.}
% 35J93, 53C50, 35J66, 35B45, 49Q05, 53A10, 78A25.\par
% \noindent\footnotesize\textbf{Key words and phrases.}
% Prescribed Lorentzian mean curvature; space-like graphs; capillary boundary
% condition; free boundary; Born--Infeld equation.
\normalsize

\section{Introduction}

\subsection*{The problem}

In this paper, we investigate existence and regularity of weak solutions to the following elliptic oblique derivative problem over a bounded Lipschitz domain $\Omega \subseteq \R^m$, $m\geq 2$,
\begin{align}\label{eq:problem_Rm}
        -\div\Bigg(\dfrac{Du}{\sqrt{1 - |Du|^2}}\Bigg)= \rho \quad \text{in $\Omega$,} \qquad -\dfrac{Du \cdot n}{\sqrt{1 - |Du|^2}} = \psi \quad \text{on $\del \Omega$},
        % \qquad \int_B u \, dx = 0,
\end{align}
where $n$ is the exterior unit normal of $\del \Omega$ and $\rho\in\mathscr M(\Omega)$ and $\psi\in\mathscr M(\del\Omega)$ are given finite signed Radon measures on $\Omega$ and $\del \Omega$ respectively. By a \emph{weak solution} to \eqref{eq:problem_Rm} in a relatively open set $\Omega' \subseteq \overline\Omega$, we mean a function $u \in \lip(\Omega)$ such that
\begin{align}
    w_u &\doteq \frac{1}{\sqrt{1 - |Du|^2}} \in L^1(\Omega'),\label{eq:energy} \\ \int_\Omega w_u Du \cdot D\eta \, dx &= \int_\Omega \eta \, d\rho - \int_{\del \Omega} \eta \, d\psi \label{eq:Rm_weak_sol} \qquad \text{for every $\eta \in \lip_c(\Omega')$},
\end{align}
where $\lip_c(\Omega')$ denotes the set of compactly supported Lipschitz functions on $\Omega'$. Here and below, test functions in $\lip_c(\Omega')$ are extended by zero to $\overline\Omega$. 

Problem \eqref{eq:problem_Rm} is relevant both in geometry and physics. From the geometrical viewpoint, the graph of a classical solution to \eqref{eq:problem_Rm} is a (strictly) space-like graph in Minkowski spacetime $\Mink$ whose mean curvature is $-\rho$ and that meets the cylinder $\R\times \del\Omega$ with hyperbolic angle $\gamma = -\operatorname{arcsinh} \psi$, see \cref{sec:preliminaries} below. We say that a function $u \in \lip(\Omega)$ is
\begin{itemize}
    \item \emph{weakly space-like} if $\norm{Du}_{L^\infty(\Omega)}\leq 1$;
    \item \emph{space-like} if, for every distinct $x,y \in \Omega$ such that the segment $\overline{xy}$ is contained in $\Omega$, it holds
    \begin{align*}
        |u(x)-u(y)|<|x-y|;
    \end{align*}
    \item \emph{strictly space-like} if $u \in C^1(\Omega)$ and $|Du|<1$ in $\Omega$;
    \item \emph{uniformly space-like} if $w_u \in L^\infty(\Omega)$.
\end{itemize}

Following the geometric interpretation, the function $w_u$ is often called the \emph{tilt function}, as it measures the inclination of the graph with respect to the vertical direction.
% and blows to $+ \infty$ as the space-likeness of the graph degenerates

On the physical side, \eqref{eq:problem_Rm} can be interpreted as an electrostatic problem within the framework of Born--Infeld's Model of Electrodynamics, a non-linear classical theory of the electromagnetic field proposed in the 30s to overcome the problem of infinite energy associated to point charges according to Maxwell's Electromagnetism. The chief feature of Born--Infeld's Electromagnetism is that the electromagnetic field cannot exceed a universal bound (normalized to $1$), much as in Special Relativity the speed of a massive particle cannot exceed the speed of light. Specializing to electrostatics, this means that the electrostatic field $\bm E = - Du$, $u$ being the electrostatic potential, is subject to the constraint $|\bm E| < 1$. The first equation in \eqref{eq:problem_Rm} is then a nonlinear analogue of Poisson’s equation. The electric displacement field is $\mathbf D=-Du/\sqrt{1-\lvert Du\rvert^2}$, and the boundary condition prescribes its normal component.
% This condition naturally arises in many physical situations, \oorange{when a superficial electrostatic charge $\psi = \sinh\gamma$ is placed on the boundary, see below. APPENDICE?}.
Physically, the function $w_u$ is closely related to the electrostatic energy density of the electric field (see for instance the Appendix in \cite{BIMM} or \cite[Chapter 3]{ManiscalcoThesis2026}) and is therefore often simply called \emph{energy density}.

In light of the physical interpretation, it is natural to investigate the problem under low regularity assumptions on $\rho$ and $\psi$, so as to model physically relevant charge distributions. The Dirichlet counterpart of \eqref{eq:problem_Rm}
\begin{align}\label{eq:Dirichlet}
        -\div\Bigg(\dfrac{Du}{\sqrt{1 - |Du|^2}}\Bigg)= \rho \quad \text{in $\Omega$,} \qquad u = \varphi \quad \text{on $\del \Omega$},
        % \qquad \int_B u \, dx = 0,
\end{align}
where $\varphi \in C(\del\Omega)$ is a given function, has attracted increasing attention in this direction. Existence and regularity were first discussed by Bartnik and Simon in their seminal work \cite{BartnikSimon1982}, where they address the case $\rho \in L^\infty$. Less regular sources were investigated by Klyachin and Miklyukov \cite{Klyachin03,KlyMik92,KlyMik93}, in which a finite superposition of Dirac deltas is studied. Later, Byeon, Ikoma, Malchiodi and Mari in \cite{BIMM} investigated the case where $\rho$ is a Radon measure. More recently, De Filippis and Mingione \cite{DeFilippisMingione26}
developed a regularity theory for variational solutions with unbounded
prescribed mean curvature, obtaining sharp results in the scale of
Lorentz spaces.

% The global problem on the whole problem has also been extensively studied, \oorange{see the contributions of ..}. See \cite{BIMM} for a complete account of literature to date.

By contrast, non-Dirichlet boundary problems for the Lorentzian mean curvature operator constitute a quite open field of research. Bereanu, Jebelean and
Mawhin studied radial solutions to \eqref{eq:problem_Rm} with continuous
source $\rho$ in \cite{Bereanu10} (see also \cite{Bereanu11,Bereanu11bis}).
In our previous work \cite{Maniscalco26}, we treated the homogeneous
Neumann problem ($\psi=0$) on convex domains for $\rho\in L^\infty(\Omega)$.
The regularity results obtained there allow one to identify critical points
of the associated nonsmooth functionals as strictly space-like solutions.
Applications include multiplicity results for problems with indefinite
weights \cite{Ziegele26} and odd nonlinearities \cite{BoscagginZiegele26}.

Before stating our main improvements, we give a quick overview of the mathematical challenges that are to be faced in studying existence and regularity for problems involving the Lorentzian mean curvature. Problem \eqref{eq:problem_Rm} arises as the (formal) Euler--Lagrange equation of the functional
\begin{align}\label{eq:action}
    I_{\rho,\psi}(u) \doteq \int_\Omega \left( 1 -  \sqrt{1 - |Du|^2} \right) \, dx - \int_\Omega u \, d\rho + \int_{\del\Omega} u \, d\psi, \qquad u \in \mathscr A(\Omega),
\end{align}
where
\begin{align}\label{eq:achronal}
    \mathscr A(\Omega) = \left\{ \ u \in \lip(\Omega) \quad | \quad  \norm{Du}_{L^\infty(\Omega)} \leq 1, \quad \int_\Omega u \, dx= 0 \ \right\}.
\end{align}
The functional $I_{\rho,\psi}$ is strictly convex and, if $|\rho|$ and $|\psi|$ are finite measures, then it is bounded below on $\mathscr A(\Omega)$. If $\rho$ and $\psi$ satisfy the compatibility condition
\begin{align}\label{eq:compatibility}
        \int_\Omega  d\rho = \int_{\del \Omega} d\psi,
\end{align}
alongside with some mild assumptions on $\rho$ and $\psi$, then $I_{\rho,\psi}$ has a unique minimiser in $\mathscr A(\Omega)$ (cf. \cref{lem:magi} below). Such minimiser is hence a natural candidate to being a weak solution to \eqref{eq:problem_Rm}. However, the lack of smoothness of the Lagrangian density $\mathscr L(p) \doteq 1-\sqrt{1 - |p|^2}$ when $|p|=1$ prevents us from concluding directly that a minimiser of $I_{\rho,\psi}$ is actually a weak solution. One road to solvability consists in establishing a priori gradient bounds $\norm{Du}_{L^\infty(\Omega)} \leq 1 - \theta$ for the minimiser $u$, like we did in \cite{Maniscalco26}. However, in most cases, such a control is not possible and one must rely on more subtle integral bounds on the tilt function $w_u$.

A major obstacle that might prevent the minimiser $u$ from being a weak solution is the possible existence of a \emph{light segment} on its graph: there might exist distinct points $x,y \in \Omega$ such that $\overline{xy}\Subset \Omega$ and 
\begin{align*}
    |u(x) - u(y)| = |x-y|.
\end{align*}
This was already acknowledged as a technical and substantial difficulty by Bartnik and Simon in \cite{BartnikSimon1982}. For $m\geq3$, Byeon, Ikoma, Malchiodi and Mari constructed a compactly supported smooth weak solution with a light segment and source $f\in L^q(\mathbb R^m)$ for every $1\leq q<m-1$ (see \cite[Theorem 1.7]{BIMM}). Adding two opposite point charges at the endpoints of that segment leaves the minimiser unchanged, but the same function then fails to solve the equation for the modified source $\rho \doteq f \, dx + (\delta_y - \delta_x)$ (where the endpoints are labelled so that $u(y)>u(x)$, see \cite[Theorem 1.9]{BIMM} and \cite[Corollary 1.10]{BIMM}). For sources of the latter type, the authors also proved some solvability criteria: if $\Sigma\Subset \Omega$ is such that $\mathscr H^1(\Sigma) = 0$ and
\begin{align}\label{eq:rho_BIMM}
    \rho = f \, dx + \nu, \qquad f \in L^1(\Omega) \cap L^2_\loc(\Omega\setminus\Sigma), \qquad \spt \nu \subseteq \Sigma
\end{align}
then the minimiser weakly solves the Dirichlet problem in $\Omega$ if $m = 2$ while, if $m\geq 3$, it is a weak solution \emph{outside the set of light segments} (cf. \cite[Theorems 1.11 and 1.14]{BIMM}). If $m\geq 4$, their counterexample has source of the form \eqref{eq:rho_BIMM}, hence the minimiser does not need to be a weak solution. The case $m = 3$ is open.

% A major obstacle to solvability, as the example constructed in \cite[Theorem 1.7]{BIMM} shows, is the possible presence of a \emph{light segment} on the minimiser $u$, that is, a segment $\overline{xy}\subset\Omega$ such that
% \begin{align*}
%     |u(x) - u(y)| = |x-y|.
% \end{align*}
% The terminology comes from the fact that, over a light segment, the graph of $u$ is a light-like (or null) hypersurface of Minkowski spacetime $\Mink$.

\subsection*{Main results}

In this paper we mainly deal with the case where $\Omega = B = B_R(o)$ is a ball of $\R^m$. We study properties of the minimiser of $I_{\rho,\psi}$ under progressively weaker assumptions on $\rho$ and $\psi$. We begin with the case of bounded $\rho$ and Lipschitz $\psi$, for which we are able to prove uniform space-likeness of the minimiser. Hereafter, we will write $w = w_u$ if there is no risk of confusion.

\begin{introthm}\label{thm:existence_lip}
    Let $\rho \in L^\infty(B)$ and $\psi \in \lip(\del B)$. Then the minimiser $u$ of $I_{\rho,\psi}$ in $\mathscr A(B)$ is a weak solution to \eqref{eq:problem_Rm} if and only if $\rho$ and $\psi$ satisfy the compatibility condition \eqref{eq:compatibility}. In this case $u$ is of class $W^{2,2}(B)\cap C^{1,\beta}(\overline B)$ for some $\beta \in (0,1)$ and enjoys the gradient bound
    \begin{align*}
        \sup_B w \leq C\left(m,R,\norm{\rho}_{L^\infty(B)}, \norm{\psi}_{\lip(\del B)}\right).
    \end{align*}
    Moreover,
    \begin{gather}\label{eq:uniform_second_control}
        \int_B \left\{ w^2|D^2u|^2 + |D\ln w|^2 + (Du \cdot Dw)^2 \right\} \, dx \leq C \left(m,R,\norm{\rho}_{L^2(B)}, \norm{D^\top \psi}_{L^2(\del B)} \right)
    \end{gather}
    and any other weak solution to \eqref{eq:problem_Rm} differs from $u$ by an additive constant. 
\end{introthm}

When treating less regular sources, we cannot expect uniform space-likeness of the minimiser to hold. The loss of regularity already occurs for homogeneous boundary data, as the following example shows.

\begin{ex}
    In the unit ball, if $m\geq3$, set
    \begin{align*}
        u(x)=r-\tfrac23r^{3/2}, \qquad r \doteq |x|.
    \end{align*}
    Then $u'(1)=0$ and $w_u=(2\sqrt r-r)^{-1/2}$. A direct computation gives $\rho=-\div(w_uDu)=O(r^{-5/4})\in L^2(B)$. Since $r^{m-1}w_u u'(r)\to0$ as $r\downarrow0$, the distributional divergence has no atom at the origin. Thus, after subtraction of its mean, $u$ is the minimiser for the data $\rho$, $\psi = 0$. Its tilt is unbounded and $u$ is not differentiable at the origin.
\end{ex}

The next example shows that uniform space-likeness can also fail for bounded boundary data, even when $\rho=0$.

\begin{ex}
    Consider the smooth, space-like function $u$ on $B_1 \subseteq \R^2$ given by
    \begin{align*}
        u(x,y) = y\tanh x
    \end{align*}
    It is clear that $|Du|=1$ at $(0,\pm1)$. A direct computation shows that $-\rho = \div (w_u Du) = 0$ in $B$ and that $\psi = - w_u Du \cdot n$ is bounded, but not continuous on $\del B$. 
\end{ex}

However, except for uniform space-likeness and $C^{1,\beta}$ regularity, the conclusions of \cref{thm:existence_lip} still hold if one only assumes $\rho \in L^2(B)$ and $\psi \in W^{1,2}(\del B)$.

\begin{introthm}\label{thm:existence_H^1}
    Let $\rho \in L^2(B)$ and $\psi \in W^{1,2}(\del B)$ satisfy the compatibility condition \eqref{eq:compatibility}. Then the minimiser $u$ of $I_{\rho,\psi}$ on $\mathscr A(B)$ is of class $W^{2,2}(B)$ and weakly solves \eqref{eq:problem_Rm}. Moreover, $u$ satisfies \eqref{eq:uniform_second_control} and
    \begin{gather}\label{eq:uniform_w_control}
        \int_B w^2 \, dx \leq C \left(m,R,\norm{\rho}_{L^2(B)}, \norm{D^\top \psi}_{L^2(\del B)} \right).
    \end{gather}
    Finally, any other weak solution to \eqref{eq:problem_Rm} differs from $u$ by an additive constant.
\end{introthm}

Our next Theorem exploits the interior regularity results of \cite{BIMM} to treat the case where $\rho$ is a measure with a singular part.

\begin{introthm}\label{thm:ultimate_thm}
    Let $\Sigma\Subset B$ be a compact set with $\Haus^1(\Sigma) = 0$ and let $\rho$ be a Radon measure on $B$ with Lebesgue decomposition of the form
    \begin{align*}
        \rho = f\, dx + \nu, \qquad f \in L^2_{\loc}(B\setminus \Sigma) \cap L^1(B), \qquad \spt \nu \subseteq \Sigma.
    \end{align*}
    Assume also that $\psi \in L^q(\del B)$ for some $q>1$ and that the compatibility condition \eqref{eq:compatibility} with $\rho$ holds. Let $u$ be the unique minimiser of $I_{\rho,\psi}$ in $\mathscr A(B)$ and let
    \begin{align*}
        K_\psi^\rho \doteq \overline{\bigcup \, \set{\, \overline{xy} \ | \ x,y \in B, \ x \neq y, \ |u(x) - u(y)| = |x - y| \, }},
    \end{align*}
    where the closure is taken in $\overline B$. Then
    \begin{enumerate}[label=(\roman*)]
        \item\label{item:ultimate_weak_solution} the minimiser $u$ of $I_{\rho,\psi}$ is of class $W_{\loc}^{2,2}(B\setminus (\Sigma\cup K_\psi^\rho))$, and is a weak solution to \eqref{eq:problem_Rm} in $\overline B\setminus K_\psi^\rho$;
        \item\label{item:ultimate_higher_regularity} if $\Omega'\Subset B\setminus (\Sigma \cup K_\psi^\rho)$ and $\rho \in L^\infty(\Omega')$, then $u$ is of class $C_{\loc}^{1,\beta}(\Omega')$ for some $\beta \in (0,1)$ and locally uniformly space-like in $\Omega'$ (i.e. $w \in L_{\loc}^\infty(\Omega')$). In particular, if $\rho \in C^\infty(\Omega')$, then so is $u$;
        \item\label{item:m=2} if moreover $m = 2$, the set $K_\psi^\rho$ is empty.
    \end{enumerate}
\end{introthm}

The conclusions of \cref{thm:ultimate_thm} can be improved if one is able to exclude the possibility for the minimiser to have light segments. In our next and last result, we are able to do so on a general bounded convex domain.

\begin{introthm}\label{thm:no_light_segment}
    Let $\Omega \subset \R^m$ be a bounded, convex domain, and let $\rho \in L^{m-1}(\Omega)$ and $\psi \in L^{m-1}(\del\Omega)$ satisfy the compatibility condition \eqref{eq:compatibility}. Then the unique minimiser of $I_{\rho,\psi}$ in $\mathscr A(\Omega)$ does not have light segments in the interior of $\Omega$, i.e. $K_\psi^\rho = \varnothing$.
\end{introthm}

Our proof builds on the variational argument of De Filippis and Mingione for rough anti-peeling \cite[Theorem~1.1]{DeFilippisMingione26}. In the capillary setting the boundary trace is free, and the additional boundary term in the variational comparison is controlled using the convexity of $\Omega$ and the integrability of the prescribed flux $\psi$.

Combining \cref{thm:no_light_segment} with our previous regularity results we obtain the following

\begin{cor*} Assume $\rho \in \mathscr M(B)$ and $\psi \in \mathscr M(\del B)$ are compatible.
\begin{enumerate}
    \item In dimension $m\geq3$, if $\rho \in L^{m-1}(B)$ and $\psi \in L^{m-1}(\del B)$, then the minimiser does not have light segments and weakly solves \eqref{eq:problem_Rm} in the whole of $\overline B$.
    \item In dimensions $m = 2,3$, under the hypotheses of \emph{\cref{thm:existence_H^1}}, the minimiser does not have light segments.
\end{enumerate}
\end{cor*}

The last conclusion fails in dimension $m\geq 4$. By the aforementioned counterexample \cite[Theorem 1.7]{BIMM}, there exists $\rho \in L^2(B)$ such that the minimiser of $I_{\rho,\psi}$, with $\psi = 0$, has a light segment.

\subsection*{Strategy} To prove that the minimiser is a weak solution we adopt the following strategy. We smooth the data by convolution and adjust the interior source by a constant so that the compatibility condition holds for every $j$. Clearly
\begin{align*}
    \lim_{j\to \infty} \int_B \eta \, d\rho_j = \int_B \eta \, d\rho, \qquad \lim_{j\to \infty} \int_{\del B} \eta \, d\psi_j = \int_{\del B} \eta \, d\psi, \qquad \text{for every $\eta \in \lip(B)$.}
\end{align*}

Our \emph{a priori} gradient estimate in \cref{thm:full_gradient_estimate} together with Lieberman's theory for oblique derivative elliptic problems ensures that for every $j$ the minimiser $u_j$ of $I_{\rho_j,\psi_j}$ is smooth and classically solves the approximating problem
\begin{align*}
    -\div\left(w_j Du_j\right)= \rho_j \quad \text{in $B$,} \qquad - w_j Du_j \cdot n = \psi_j \quad \text{on $\del B$}.
        % \qquad \int_B u \, dx = 0,
\end{align*}

As the $u_j$ converge to $u$ (cf. \cref{lem:magi}), the latter is a weak solution to \eqref{eq:problem_Rm} in a relatively open $\Omega \subseteq \overline B$ if and only if
\begin{align}\label{eq:pass_to_the_limit}
    \lim_{j \to \infty} \int_{B} w_j Du_j \cdot D\eta \, dx = \int_B wDu \cdot D\eta \, dx, \qquad \text{for every $\eta \in \lip_c(\Omega)$.}
\end{align}
By Vitali's convergence Theorem \cite[Theorem 3.1.9]{Willem2013}, \eqref{eq:pass_to_the_limit} follows provided that $\set{w_j}$ is locally uniformly integrable in $\Omega$ (when $\Omega = \overline B$, uniform integrability is also necessary, as follows by testing the approximating equations with $u_j$ and the limit equation with $u$).
% Recall that this means that, for every $\eps>0$ there exists $\alpha>0$ such that for every measurable set $A \subseteq B$ it holds
%     \begin{align*}
%         |A|<\alpha \qquad \Longrightarrow \qquad \sup_j \int_A w_j \, dx < \eps.
%     \end{align*}
In turn, by de la Vallée-Poussin's Theorem \cite[Theorem~3.1.10]{Willem2013}, this is guaranteed by the existence of a convex, increasing function $h:[0,\infty) \to [0,\infty)$ such that
\begin{align}\label{eq:poussin}
    \lim_{t\to\infty} \frac{h(t)}{t} = \infty, \qquad \sup_j \int_{\Omega'} h(w_j) \, dx < \infty, \qquad \text{for every $\Omega'\Subset \Omega$.}
\end{align}

Under the assumptions of \cref{thm:existence_H^1}, a global integral estimate for the second fundamental form of the graph of a smooth, strictly space-like function, \cref{lem:L^2_w_estimates}, applied along the sequence $u_j$, gives \eqref{eq:poussin} with $h(t) = t^2$. More precisely, it holds
\begin{align*}
    \int_B w_j^2 \, dx \leq C(m,R,\norm{\rho}_{L^2(B)},\norm{D^\top \psi}_{L^2(\del B)}).
\end{align*}
This is the key to the proof of \cref{thm:existence_H^1}.

Under the assumptions of \cref{thm:ultimate_thm}, a uniform $L^2$ bound on $D^\top \psi_j$ is no longer available. We instead apply \cite[Theorem 4.13]{BIMM}, which roughly states that, if $\rho$ is as in \cref{thm:ultimate_thm}, then
\begin{align*}
    \sup_j\int_{\Omega'} w_j \ln (1 + w_j) \, dx < \infty \qquad \text{for every $\Omega'\Subset B\setminus (\Sigma\cup K_\psi^\rho)$}
\end{align*}
thus establishing \eqref{eq:poussin} with $h(t) = t \ln (1 + t)$ and $\Omega = B\setminus (\Sigma\cup K_\psi^\rho)$. As a consequence, $\set{w_j}$ is locally uniformly integrable on $B\setminus (\Sigma\cup K_\psi^\rho)$, and thus, by the Removable Singularities result \cite[Theorem 4.2]{BIMM}, $\set{w_j}$ is locally uniformly integrable in the whole of $B\setminus K_\psi^\rho$.

In order to conclude that \eqref{eq:pass_to_the_limit} holds, we only need to control the integral of $w_j$ near the boundary, uniformly in $j$. This is achieved in \cref{thm:boundary_energy_estimate}, which gives
\begin{align*}
    \int_{U_\delta} w_j \, dx \leq C\left(\delta^{\alpha_q} + \int_{U_{2\delta}} |\rho_j| \, dx \right), \qquad C = C(m,R,q,\norm{\rho_j}_{L^1(B)},\norm{\psi_j}_{L^q(\del B)}),
\end{align*}
where $U_\delta = B_R \setminus \overline{B_{R - \delta}}$ and $\alpha_q>0$ only depends on $q$. This immediately yields local uniform integrability of $\set{w_j}$ on the whole of $\overline B\setminus K_\psi^\rho$, thus proving \cref{thm:ultimate_thm}.

The gradient estimate \cref{thm:full_gradient_estimate} and the integral control of the tilt near the boundary, \cref{thm:boundary_energy_estimate}, constitute the technical core of the paper. Both use the rotational symmetries of the ball. However, we believe that many of the techniques and results we present could be adapted to study problem \eqref{eq:problem_Rm} on $C^{1,1}$-regular, convex, bounded domains.

\subsection*{Acknowledgements} AI tools were used for double checking proofs and computations. When asked to verify the proof of \cref{thm:boundary_energy_estimate} below, an AI found and fixed a gap. As a consequence, Step 3 and 4 in the proof of \cref{thm:boundary_energy_estimate} (and the Appendix, which is closely related) have been written with the assistance of the machine.

\section{Preliminaries from Lorentzian Geometry}\label{sec:preliminaries}

% This first Section explains the geometric and physical interpretations of \eqref{eq:problem_Rm} and fixes the conventions that will be followed throughout the paper.
Hereafter, derivatives of a smooth function $f$ on $\R^m$ with respect to any fixed coordinate system $\set{x^i}$ will be denoted by
\begin{align*}
    f_i \doteq \frac{\del f}{\del x^i}, \qquad f_{ij} \doteq \frac{\del^2 f}{\del x^i \del x^j}, \qquad \text{and so on.}
\end{align*}
Einstein's summation convention will be used throughout the paper. Indices will always be raised with respect to the Euclidean metric $e = e_{ij} \, dx^i \otimes dx^j$, so we will always write $f^i \doteq e^{ij}f_j$, where $e^{ij}$ is the inverse matrix of $e_{ij}$. We will usually use Greek and Latin indices with the following ranges
\begin{align*}
    \mu, \nu, \lambda,\dots \in \set{0,1,2,\dots,m}, \qquad i,j,k,\dots \in \set{1,2,\dots,m}, \qquad a,b,c,\dots \in\set{2,\dots,m}.
\end{align*}

Minkowski spacetime $\Mink$ is $\R\times \R^m$ endowed with the pseudo-scalar product
\begin{align*}
    \langle \cdot, \cdot \rangle = -dx^0 \otimes dx^0 + dx^1 \otimes dx^1 + \dots + dx^m \otimes dx^m
\end{align*}
where $\set{dx^\mu}$ is the dual basis of the standard basis $\set{\del_\mu}$ of $\R\times \R^m$. We will denote by $\babla$ the Levi-Civita connection of $\Mink$. Throughout this section, let $\Omega$ be a convex domain with smooth boundary and let $u\in C^2(\overline\Omega)$ be strictly space-like in $\overline\Omega$. Consider its graph in $\Mink$:
\begin{align*}
    F: \overline\Omega \longrightarrow \Mink, \qquad F(x) \doteq (u(x),x), \qquad M \doteq F(\overline\Omega).
\end{align*}
The future-pointing, unit normal $N$ to the graph of $u$ is well defined and time-like. In terms of $u$, it is given by
\begin{align*}
    N \doteq w(Du + \del_0), \qquad w = \frac{1}{\sqrt{1-|Du|^2}} = - \langle N,\del_0\rangle.
\end{align*}
As anticipated, $w$ is the hyperbolic cosine of the hyperbolic angle between $N$ and $\del_0$, hence the name \emph{tilt function.} Call $g \doteq F^*\langle\cdot,\cdot\rangle$ the \emph{graph metric} of $u$, so that $(\overline{\Omega},g)$ is isometric to $M$ with the induced metric. Since the differential of the graph map $F$ satisfies $F_*\del_i = \del_i + u_i \del_0$ for every coordinate system $\set{x^i}$, we see that $g$ and its inverse $g^{-1}$ have components
\begin{align*}
    g_{ij} = e_{ij} - u_iu_j, \qquad g^{ij} = e^{ij} + w^2 u^i u^j,
\end{align*}
where $e_{ij}$ and $e^{ij}$ are the Euclidean metric and its inverse in the chosen coordinate system. We denote by $D$ the Euclidean Levi-Civita connection on $\R^m$ and $\nabla$ that of $(\overline\Omega,g)$. The Euclidean and graph gradients of a smooth function $f$ are defined by
\begin{align*}
    D f \doteq e^{ij}f_j \del_i, \qquad \nabla f \doteq g^{ij} f_j \del_i.
\end{align*}
We will always compute the norm of a gradient with respect to the same metric it refers to, and we will simply write
\begin{align*}
    |Df|^2 \doteq e^{ij}f_i f_j, \qquad |\nabla f|^2 = g^{ij} f_i f_j.
\end{align*}
By computing the determinant of $g_{ij}$ one can verify that $w$ is the volume density of the Lebesgue measure on $\overline\Omega$ with respect to the Riemannian volume element $d\mu$ of $g$:
\begin{align}\label{eq:volumes}
    dx = w \, d\mu.
\end{align}

The Levi-Civita connection of the graph and that of the ambient $\Mink$ are related by the Gauss formula
\begin{align}\label{eq:gauss_equation}
    \babla_{F_*X}F_*Y = F_*(\nabla_XY) + \II(X,Y)N, \qquad \text{for all $X,Y \in\X(\overline\Omega)$}
\end{align}
where $\II$ is the \emph{second fundamental form} of $F$, the $(0,2)$-tensor on $\overline\Omega$ defined by
\begin{align*}
    \II(X,Y) = - \langle \babla_{F_*X}F_*Y,N\rangle = \langle F_*Y, \babla_{F_*X}N \rangle.
\end{align*}
From its definition we obtain $\II = w D^2u$, where $D^2u$ is the Euclidean Hessian of $u$. The metric equivalent $(1,1)$ tensor $A$ associated to $\II$ is called the \emph{shape operator} and is defined by
\begin{align*}
    \II(X,Y) = g(AX,Y), \qquad \text{for every $X,Y \in\X(\overline\Omega)$}.
\end{align*}
The trace of the shape operator is the (unnormalised) \emph{mean curvature} in direction $N$, and is denoted by $H \doteq \tr A = \tr_g \II = g^{ij}\II_{ij}$. In terms of $u$,
\begin{align*}
    H = \div\left(wDu\right).
\end{align*}
In intrinsic terms, 
\begin{align}\label{eq:mean_curv_intrinsic}
    \Delta_g F = HN, \qquad \Delta_g u = Hw, \qquad \Delta_g x^i = Hwu^i,
\end{align}
where $\Delta_g$ is the Laplace--Beltrami operator of $g$ and $\set{x^i}$ are the standard Cartesian coordinates.
% The mean curvature can equivalently be defined as the tangential divergence $H = \div_\Sigma N$ of the normal vector field, where
% \begin{align*}
%     \div_\Sigma X \doteq g^{ij}\langle \babla_{F_*{\del_j}} X, F_*\del_j \rangle, \qquad X \in \X(\Mink).
% \end{align*}

The boundary $\del M$ can be treated the same way by considering the graph of $u|_{\del\Omega}$ in the product spacetime $\R\times \del\Omega$ instead of in $\Mink$. Let $n$ be the outward unit normal of $\del\Omega$. Choose adapted coordinates $\set{x^i}$ in $\Omega$ so that $x^1 \equiv 0$ on $\del \Omega$ and $\del_1 = n$. Then the metric of $\del \Omega$ has components $e_{ab}$. Denote by $D^\top$ the Levi-Civita connection of $\del\Omega$, so that, in particular,
\begin{align*}
    D^\top f \doteq e^{ab} f_a \del_b, \qquad f \in C^1(\del \Omega).
\end{align*}
For a function $f \in C^1(\overline \Omega)$ it holds
\begin{align*}
    Df = D^\top f + (Df \cdot n)n \qquad \text{along $\del \Omega$.}
\end{align*}
The graph metric of $\del M$ and its inverse are
\begin{align}\label{eq:boundary_graph_metric}
    g_{ab} = e_{ab} - u_a u_b, \qquad g^{ab} = e^{ab} + w_\top^2 (D^\top u)^a (D^\top u)^b, \qquad w_\top \doteq \frac{1}{\sqrt{1 - |D^\top u|^2}}.
\end{align}
As above, if we let $d\sigma$ denote the Riemannian volume element of $\del M$, then $w_\top$ is the density of the surface measure $d\omega$ on $\del \Omega$ with respect to $d\sigma$:
\begin{align}\label{eq:boundary_volumes}
    d\omega = w_\top \, d\sigma.
\end{align}

Letting $\bm n = (0,n)$ be the exterior unit normal of the cylinder $\R \times \Omega \subset \Mink$, we let $\gamma$ be the hyperbolic angle between $M$ and $\R \times \del\Omega$, that is $\langle N,\bm n\rangle \doteq \sinh \gamma$, or, in terms of $u$,
\begin{align*}
    \sinh \gamma = wDu \cdot n.
\end{align*}
Observe that this implies
\begin{align*}
    w = w_\top \cosh \gamma.
\end{align*}
Let now $\nu$ be the exterior conormal vector of $\del M$ inside $M \cong (\overline\Omega,g)$. A computation and the capillary condition imply
\begin{align}\label{eq:conormal}
    \nu = \frac{n + w\sinh\gamma \,  Du}{\cosh\gamma}.
\end{align}
Suppose now that $\Omega=B_R$, and let $r$ denote the distance from its centre. One can verify that 
\begin{align}
    \langle \nabla u,\nu\rangle = w\tanh\gamma, \qquad \langle \nabla r, \nu\rangle = \cosh \gamma
\end{align}
and
\begin{align}\label{eq:useful_identities}
    |\nabla u|^2 = w^2 - 1, \qquad w^2(\del_ru)^2 = |\nabla r|^2 - 1, \qquad \langle\nabla u,\nabla r\rangle = w^2 \del_r u.
\end{align}
Combining these identities, we deduce
\begin{align}\label{eq:intrinsic_normal_derivatives}
    \langle \nabla u,\nu\rangle \, d\sigma = \sinh\gamma \, d\omega, \qquad \langle \nabla r,\nu\rangle \, d\sigma = \frac{\cosh \gamma}{w_\top} \, d\omega.
\end{align}

\section{Absence of light segments} 

This section is devoted to proving that, if $\rho \in L^{m-1}(\Omega)$ and $\psi\in L^{m-1}(\del\Omega)$ are compatible, then the Born--Infeld action $I_{\rho,\psi}$ defined in \eqref{eq:action} has a unique minimiser in the set $\mathscr A(\Omega)$, and that this minimiser has no interior light segments. We first recall some standard facts about minimisers of $I_{\rho,\psi}$.

\begin{lem}\label{lem:magi}
    Let $\Omega \subset \mathbb R^m$ be a bounded domain with Lipschitz boundary.
    Let $\rho,\rho_j \in \mathscr M(\Omega)$ and
    $\psi,\psi_j \in \mathscr M(\del\Omega)$ be finite signed Radon
    measures such that $(\rho,\psi)$ and every $(\rho_j,\psi_j)$
    satisfy the compatibility condition \eqref{eq:compatibility}.
    Then:
    \begin{enumerate}[label=(\roman*),ref=\roman*]
        \item\label{eq:melchior}
        there exists a unique minimiser $u$ of $I_{\rho,\psi}$
        in $\mathscr A(\Omega)$;
        \item\label{eq:balthasar}
        if $u \in \mathscr A(\Omega)$ is a weak solution to
        \eqref{eq:problem_Rm}, then it is the minimiser of
        $I_{\rho,\psi}$ in $\mathscr A(\Omega)$;
        \item\label{eq:casper}
        if
        \[
            \rho_j \overset{*}{\rightharpoonup} \rho
            \quad\text{in }\mathscr M(\overline\Omega),
            \qquad
            \psi_j \overset{*}{\rightharpoonup} \psi
            \quad\text{in }\mathscr M(\del\Omega),
        \]
        where the measures on $\Omega$ are extended by zero to
        $\overline\Omega$, then the minimisers $u_j$ of
        $I_{\rho_j,\psi_j}$ converge to $u$ strongly in
        $C(\overline\Omega)\cap W^{1,p}(\Omega)$ for every
        $1 \leq p < \infty$, and weakly-* in $\lip(\Omega)$.
    \end{enumerate}
\end{lem}

\begin{proof}

    Set
    \begin{align}\label{eq:J}
        J(v) \doteq \int_\Omega \mathscr L(Dv)\,dx,
        \qquad
        \mathscr L(p) \doteq 1-\sqrt{1-|p|^2},
        \quad |p|\leq 1,
    \end{align}
    and write
    \[
        \ell_{\rho,\psi}(v)
        \doteq \int_\Omega v\,d\rho
        -\int_{\del\Omega}v\,d\psi,
        \qquad
        I_{\rho,\psi}(v)=J(v)-\ell_{\rho,\psi}(v).
    \]

    Since $\Omega$ is Lipschitz, every $v\in\mathscr A(\Omega)$
    has a $1$-Lipschitz representative on $\overline\Omega$.
    Together with the zero-mean condition, this gives
    \begin{equation}\label{eq:stupid_L^infty_bound}
        \norm{v}_{L^\infty(\Omega)}
        \leq \operatorname{diam}\Omega,
        \qquad v\in\mathscr A(\Omega).
    \end{equation}
    Thus $\mathscr A(\Omega)$ is compact in $C(\overline\Omega)$
    by the Arzel\`a--Ascoli Theorem. Moreover, $0\leq J(v)\leq|\Omega|$
    and
    \[
        |\ell_{\rho,\psi}(v)|
        \leq \operatorname{diam}\Omega
        \bigl(|\rho|(\Omega)+|\psi|(\del\Omega)\bigr),
    \]
    so $I_{\rho,\psi}$ is bounded on $\mathscr A(\Omega)$.

    Let $\{v_k\}$ be a minimising sequence. Up to a subsequence,
    \[
        v_k\to u \quad\text{in }C(\overline\Omega),
        \qquad
        v_k\overset{*}{\rightharpoonup}u
        \quad\text{in }\lip(\Omega),
    \]
    for some $u\in\mathscr A(\Omega)$.
    Since the measures are finite,
    $\ell_{\rho,\psi}(v_k)\to\ell_{\rho,\psi}(u)$.
    By convexity of $\mathscr L$, the functional $J$ is weakly
    lower semicontinuous in $W^{1,2}(\Omega)$ on
    $\mathscr A(\Omega)$. Consequently,
    \[
        I_{\rho,\psi}(u)
        \leq \liminf_{k\to\infty}I_{\rho,\psi}(v_k)
        =\inf_{\mathscr A(\Omega)}I_{\rho,\psi}.
    \]
    This proves existence. Strict convexity of $\mathscr L$
    implies that any two minimisers have the same gradient;
    since $\Omega$ is connected and both have zero mean, they
    coincide. This proves \eqref{eq:melchior}.

    Suppose now that $u\in\mathscr A(\Omega)$ is a weak solution.
    For every $v\in\mathscr A(\Omega)$, testing
    \eqref{eq:Rm_weak_sol} with $v-u$ and using convexity gives
    \begin{align*}
        \ell_{\rho,\psi}(v-u)
        &=
        \int_\Omega D\mathscr L(Du)\cdot(Dv-Du)\,dx\\
        &\leq
        \int_\Omega
        \bigl(\mathscr L(Dv)-\mathscr L(Du)\bigr)\,dx
        =J(v)-J(u).
    \end{align*}
    Hence $I_{\rho,\psi}(u)\leq I_{\rho,\psi}(v)$, proving
    \eqref{eq:balthasar}.

    Finally, write
    \[
        I_j\doteq I_{\rho_j,\psi_j},
        \qquad I\doteq I_{\rho,\psi},
        \qquad
        \ell_j\doteq\ell_{\rho_j,\psi_j},
        \qquad \ell\doteq\ell_{\rho,\psi}.
    \]
    Weak-* convergence and the uniform boundedness principle yield
    \[
        M\doteq\sup_j
        \bigl(|\rho_j|(\Omega)+|\psi_j|(\del\Omega)\bigr)<\infty.
    \]
    In particular, if $z_j\to z$ in $C(\overline\Omega)$, then
    \begin{align}
        |\ell_j(z_j)-\ell(z)|
        &\leq
        M\norm{z_j-z}_{C(\overline\Omega)}
        +|\ell_j(z)-\ell(z)|
        \longrightarrow 0.
        \label{eq:linear_terms_convergence}
    \end{align}

    By compactness, up to a subsequence,
    $u_j\to v$ in $C(\overline\Omega)$ and
    $u_j\overset{*}{\rightharpoonup}v$ in $\lip(\Omega)$
    for some $v\in\mathscr A(\Omega)$.
    Lower semicontinuity, \eqref{eq:linear_terms_convergence},
    and the minimality of $u_j$ give
    \[
        I(v)
        \leq\liminf_{j\to\infty}I_j(u_j)
        \leq\limsup_{j\to\infty}I_j(u_j)
        \leq\lim_{j\to\infty}I_j(u)
        =I(u).
    \]
    By uniqueness, $v=u$. Since the same argument applies to
    every subsequence, the whole sequence converges to $u$
    in $C(\overline\Omega)$ and weakly-* in $\lip(\Omega)$.
    Moreover,
    \[
        I_j(u_j)\longrightarrow I(u),
        \qquad
        I(u_j)
        =I_j(u_j)+\ell_j(u_j)-\ell(u_j)
        \longrightarrow I(u).
    \]

    To obtain strong convergence of the gradients in every Lebesgue space (except $L^\infty(\Omega)$), observe that
    \[
        D^2\mathscr L(p)
        =
        \frac{\id}{\sqrt{1-|p|^2}}
        +\frac{p\otimes p}{(1-|p|^2)^{3/2}}
        \geq\id,
        \qquad |p|<1.
    \]
    Thus $p\mapsto\mathscr L(p)-|p|^2/2$ is convex on the open
    unit ball and, by continuity, also on its closure. In particular,
    \[
        \mathscr L\left(\frac{p+q}{2}\right)
        \leq
        \frac{\mathscr L(p)+\mathscr L(q)}2
        -\frac18|p-q|^2,
        \qquad |p|,|q|\leq1.
    \]
    Since $(u_j+u)/2\in\mathscr A(\Omega)$ and the lower-order
    terms are linear, the minimality of $u$ gives
    \[
        I(u)
        \leq I\left(\frac{u_j+u}{2}\right)
        \leq
        \frac{I(u_j)+I(u)}2
        -\frac18\int_\Omega|Du_j-Du|^2\,dx.
    \]
    Therefore,
    \[
        \int_\Omega|Du_j-Du|^2\,dx
        \leq4\bigl(I(u_j)-I(u)\bigr)
        \longrightarrow0.
    \]
    For $2\leq p<\infty$, the bound $|Du_j-Du|\leq2$ yields
    \[
        \int_\Omega|Du_j-Du|^p\,dx
        \leq2^{p-2}\int_\Omega|Du_j-Du|^2\,dx
        \longrightarrow0.
    \]
    For $1\leq p<2$, the conclusion follows from H\"older's
    inequality. Together with uniform convergence, this proves
    \eqref{eq:casper}.
\end{proof}

Now we move to the no-light-segment Theorem. The proof relies on a combination of ideas in \cite{Maniscalco26} with the proof of the rough anti-peeling Theorem of De Filippis and Mingione \cite[Theorem 1.1]{DeFilippisMingione26}.

\begin{proof}[Proof of \cref{thm:no_light_segment}]
    Call $u$ the minimiser and assume by contradiction that it has a light segment, i.e. there exist $o \in \Omega$, $p \in \S^{m-1}$ and $l > 0$ such that
    \begin{align*}
        u(o - tp) = u(o) - t, \qquad \text{for all $t \in [0,l]$,} \qquad S \doteq \set{o-tp \ | \ t \in [0,l]}\Subset \Omega.
    \end{align*}

    \textit{Step 1. Construction of CMC barriers.} Let
    \begin{align*}
        d \doteq \dist(o,\del\Omega), \qquad r\doteq |x-o|, \qquad R \doteq \max_{\overline\Omega}r, \qquad a_0 \doteq 1 + \frac{R^m}{m}.
    \end{align*}
    Following \cite{BartnikSimon1982} and \cite{DeFilippisMingione26}, for every (small) $\lambda > 0$ we consider the barrier
    \begin{align*}
        b_\lambda(x) = u(o) - \int_0^{|x-o|} \frac{a(s) \, ds}{\sqrt{\lambda^2 s^{2(m-1)}+a(s)^2}}, \qquad a(s) \doteq a_0 - \frac{s^m}{m}.
    \end{align*}
    Notice that $b_\lambda$ is smooth and strictly space-like away from $o$. We compute
    \begin{align*}
        Db_\lambda = - \frac{a(r)}{\sqrt{\lambda^2 r^{2(m-1)}+a(r)^2}} \frac{x-o}{r}, \qquad \text{in $\Omega\setminus\set o$}
    \end{align*}
    and therefore
    \begin{align*}
        F_\lambda \doteq \frac{Db_\lambda}{\sqrt{1 - |Db_\lambda|^2}} = - \frac{a(r)}{\lambda r^m}(x-o).
    \end{align*}
    Since $a'(r) = - r^{m-1}$, the barriers have constant classical mean curvature outside $o$:
    \begin{align*}
        \div F_\lambda = r^{1-m} \frac{d}{dr} \left( - \frac{a(r)}{\lambda} \right) = \frac1\lambda, \qquad \text{in $\Omega\setminus\set o$}.
    \end{align*}
    Moreover, by convexity of $\Omega$, for $\Haus^{m-1}$-a.e. $z \in \del\Omega$ it holds $(z-o)\cdot n(z) \geq d$, and since $a\geq 1$ we have
    \begin{align}\label{eq:barrier_boundary_flux}
        F_\lambda(z) \cdot n(z) = - \frac{a(|z-o|)}{\lambda|z-o|^m}(z-o) \cdot n(z) \leq - \frac{k}{\lambda}, \qquad k \doteq \frac{d}{R^m}. 
    \end{align}

    \textit{Step 2. Quantitative estimates around the light segment.} Denote by $\Gamma(x) \doteq u(o) - |x-o|$ the past-pointing light cone centred at $(u(o),o)$. Notice that
    \begin{align*}
        b_\lambda(x) - \Gamma(x) = \int_0^{|x-o|} \frac{\lambda^2s^{2(m-1)} \, ds}{\sqrt{\lambda^2 s^{2(m-1)} + a(s)^2}\left(\sqrt{\lambda^2 s^{2(m-1)} + a(s)^2} + a(s) \right)}, \qquad x \in \overline\Omega.
    \end{align*}
    If we restrict ourselves to $\lambda \leq \tfrac{1}{R^{m-1}}$, which we may do from now on, then $1\leq a \leq a_0$ implies
    \begin{align}\label{eq:b_lambda_comparison}
        c_0\lambda^2 r^{2m-1} \leq b_\lambda - \Gamma \leq C_0\lambda^2 r^{2m-1}, \qquad \text{on $\overline\Omega$.}
    \end{align}
    On the other hand, by weak space-likeness
    \begin{align}\label{eq:u_comparison}
        u \geq \Gamma \qquad \text{in $\overline\Omega$}.
    \end{align}
    We introduce the function
    \begin{align*}
        v_\lambda \doteq (b_\lambda - u)_+
    \end{align*}
    that we shall use to construct a competitor for $I_{\rho,\psi}$ to achieve a contradiction. Combining \eqref{eq:u_comparison} and \eqref{eq:b_lambda_comparison} we have
    \begin{align}\label{eq:bound_on_v_lambda}
        v_\lambda \leq  C\lambda^2
    \end{align}
    for some $C = C(C_0,R)$. To achieve an estimate from below, we consider the cylinder
    \begin{align*}
        Q_\lambda \doteq \set{ \ o - tp + y \quad | \quad t \in \left[\tfrac l3, \tfrac{2l}{3}\right], \quad y \perp p, \quad |y|<\delta \lambda \ }
    \end{align*}
    where $\delta>0$ will be specified later. If $x = o - tp + y \in Q_\lambda$, weak space-likeness and $u(o-lp) = u(o) - l$ imply
    \begin{align*}
        u(o - tp + y) - \Gamma(o - tp + y) &= u(o - tp + y) - u(o - lp) + |y - tp| - l \\
        &\leq |y + (l - t)p| + |y - tp| - l \\
        &= \sqrt{|y|^2 + (l - t)^2} + \sqrt{|y|^2 + t^2} - l \\
        &\leq \left(\frac{1}{2(l-t)} + \frac{1}{2t} \right) |y|^2
    \end{align*}
    where we used the inequality $\sqrt{\alpha^2 + \beta^2} \leq \tfrac{\alpha^2}{2\beta} + \beta$. Since $\tfrac{l}{3}\leq t \leq \frac{2l}{3}$, we conclude
    \begin{align*}
        u(o - tp + y) - \Gamma(o - tp + y) \leq \frac{3}{l}|y|^2, \qquad o - tp + y \in Q_\lambda.
    \end{align*}
    Since in $Q_\lambda$ we have $r \geq \tfrac{l}{3}$, \eqref{eq:b_lambda_comparison} implies
    \begin{align*}
        b_\lambda - \Gamma \geq c_1\lambda^2, \qquad \text{on $Q_\lambda$}
    \end{align*}
    where $c_1 = c_1(c_0,l)>0$ is independent of $\lambda$. Choose $\delta>0$ so that $\frac{3\delta^2}{l} \leq \tfrac{c_1}{2}$ and $\lambda\ll1$ so that $Q_\lambda \Subset \Omega$. Then
    \begin{align*}
        b_\lambda - u = b_\lambda - \Gamma + \Gamma - u \geq c_1\lambda^2 - \frac{3}{l}|y|^2 \geq \left(c_1 - \tfrac{3\delta^2}{l} \right)\lambda^2 \geq \tfrac{c_1}{2}\lambda^2 \qquad \text{on $Q_\lambda$}.
    \end{align*}
    Since the volume of the cylinder is $\tfrac{l}{3}\omega_{m-1}\delta^{m-1}\lambda^{m-1}$, where $\omega_{m-1}$ is the volume of the unit ball in $\R^{m-1}$, we obtain
    \begin{align}\label{eq:integral_bounds_v_below}
        \int_\Omega v_\lambda \, dx \geq c_2\lambda^{m+1} \qquad c_2 \doteq \tfrac{c_1\omega_{m-1}}{6}l\delta^{m-1}>0.
    \end{align}

    \textit{Step 3. Variational comparison and conclusion.} For every $\lambda,\tau>0$ consider the function
    \begin{align*}
        v_{\lambda,\tau} \doteq (b_\lambda - u - \tau)_+
    \end{align*}
    which, since $b_\lambda(o) = u(o)$, vanishes in a neighbourhood of $o$. Moreover
    \begin{align*}
        u + v_{\lambda,\tau} = \max\set{u,b_\lambda - \tau}
    \end{align*}
    is weakly space-like and, up to subtracting its mean, it belongs to $\mathscr A(\Omega)$. Since $\rho$ and $\psi$ satisfy \eqref{eq:compatibility}, adding constants does not alter the value of $I_{\rho,\psi}$, hence we can write
    \begin{align*}
        0 \geq I_{\rho,\psi}(u) - I_{\rho,\psi}(u + v_{\lambda,\tau}),
    \end{align*}
    by minimality of $u$. Notice that
    \begin{align*}
        Du + Dv_{\lambda,\tau} =
        \begin{cases}
            Db_\lambda &\qquad \text{a.e. in }\set{v_{\lambda,\tau} > 0}, \\
            Du &\qquad \text{a.e. in }\set{v_{\lambda,\tau} = 0}
        \end{cases}
    \end{align*}
    hence, since $\mathscr L$ is convex and $b_\lambda$ is strictly space-like outside $o$, we have 
    \begin{align}\nonumber
        \mathscr L(Du) - \mathscr L(Du + Dv_{\lambda,\tau}) &= \mathscr L(Du) - \mathscr L(Db_\lambda) \\ \nonumber
        &\geq - D\mathscr L (Db_\lambda) \cdot D v_{\lambda,\tau} \\
        &= - \frac{Db_\lambda \cdot Dv_{\lambda,\tau}}{\sqrt{1 - |Db_\lambda|^2}} \label{eq:convexity!}
    \end{align}
    a.e. in $\set{v_{\lambda,\tau}>0}$. On the other hand, on the complementary set $Dv_{\lambda,\tau} = 0$ a.e. and both sides vanish, hence \eqref{eq:convexity!} holds a.e. on $\Omega$ and we can write
    \begin{align*}
        0 &\geq I_{\rho,\psi}(u) - I_{\rho,\psi}(u + v_{\lambda,\tau}) \\
        &\geq - \int_\Omega F_\lambda \cdot Dv_{\lambda,\tau} \, dx + \int_\Omega v_{\lambda,\tau} \rho \, dx - \int_{\del\Omega} v_{\lambda,\tau} \psi \, d\Haus^{m-1}.
    \end{align*}
    Now take $\eps>0$ so small that $B_\eps(o) \Subset \Omega$ and $v_{\lambda,\tau}$ vanishes on $\overline{B_\eps(o)}$ and apply the divergence Theorem to obtain
    \begin{align*}
        0 \geq \int_\Omega \left( \tfrac{1}{\lambda} + \rho \right) v_{\lambda,\tau} \, dx - \int_{\del\Omega} \left(\psi + F_\lambda \cdot n\right)  v_{\lambda,\tau} \, d\Haus^{m-1}.
    \end{align*}
    By \eqref{eq:barrier_boundary_flux} this implies
    \begin{align*}
        \frac{1}{\lambda} \left( \int_\Omega v_{\lambda,\tau} \, dx + k \int_{\del\Omega} v_{\lambda,\tau}  \, d\Haus^{m-1} \right) \leq \int_{\Omega} |\rho| v_{\lambda,\tau} \, dx + \int_{\del\Omega} |\psi|v_{\lambda,\tau} \, d\Haus^{m-1}.
    \end{align*}
    As $\tau \searrow 0$, the functions $v_{\lambda,\tau}$ converge to $v_\lambda$ on $\overline\Omega$, hence, by monotone convergence
    \begin{align}\label{eq:integral_estimate}
        \frac1\lambda \left( \int_\Omega v_{\lambda} \, dx + k \int_{\del\Omega} v_{\lambda}  \, d\Haus^{m-1} \right) \leq \int_{\Omega} |\rho| v_{\lambda} \, dx + \int_{\del\Omega} |\psi|v_{\lambda} \, d\Haus^{m-1}.
    \end{align}
    Letting 
    \begin{align*}
        E_\lambda \doteq \left\{ |\rho| \geq \tfrac{1}{2\lambda} \right\}\subseteq \Omega, \qquad G_\lambda \doteq \left\{ |\psi| \geq \tfrac{k}{2\lambda} \right\} \subseteq \del\Omega,
    \end{align*}
    the integrals over $E_\lambda^c$ and $G_\lambda^c$ on the right-hand side of \eqref{eq:integral_estimate} are absorbed on the left-hand side, yielding
    \begin{align*}
        \frac{1}{2\lambda} \left( \int_\Omega v_{\lambda} \, dx + k \int_{\del\Omega} v_{\lambda}  \, d\Haus^{m-1} \right) &\leq \int_{E_\lambda} |\rho| v_{\lambda} \, dx + \int_{G_\lambda} |\psi|v_{\lambda} \, d\Haus^{m-1} \\ &\overset{\eqref{eq:bound_on_v_lambda}}{\leq} C\lambda^2 \left( \int_{E_\lambda} |\rho| \, dx + \int_{G_\lambda} |\psi| \, d\Haus^{m-1} \right) \\
        &\leq C\lambda^2 \left((2\lambda)^{m-2} \int_{E_\lambda} |\rho|^{m-1} \, dx + \left(\tfrac{2\lambda}{k}\right)^{m-2}\int_{G_\lambda} |\psi|^{m-1} \, d\Haus^{m-1} \right) \\
        &\leq C\lambda^m \left( \norm{\rho}_{L^{m-1}(E_\lambda)}^{m-1} + k^{2-m}\norm{\psi}_{L^{m-1}(G_\lambda)}^{m-1} \right).
    \end{align*}
    Combining this inequality with \eqref{eq:integral_bounds_v_below} and cancelling $\lambda^m$, we obtain
    \begin{align*}
        c_2 \leq 2C \left( \norm{\rho}_{L^{m-1}(E_\lambda)}^{m-1} + k^{2-m} \norm{\psi}_{L^{m-1}(G_\lambda)}^{m-1} \right).
    \end{align*}
    However, since $\rho\in L^{m-1}(\Omega)$ and $\psi\in L^{m-1}(\del\Omega)$, the right-hand side tends to zero as $\lambda\to 0$, contradicting \eqref{eq:integral_bounds_v_below}.
\end{proof}

\section{Gradient estimates}\label{sec:grad_est}

This section is devoted to proving our main \emph{a priori} gradient estimate:

\begin{thm}\label{thm:full_gradient_estimate}
    Let $m \geq 2$, $R> 0$ and take a zero-mean function $u \in C^3(\overline{B})$ such that $|Du| < 1$ on $\overline{B}$. Let
    \begin{align*}
        - \div(wDu) \doteq \rho \quad \text{in $B$}, \qquad - w Du \cdot n \doteq \psi \quad \text{on $\del B$}.
    \end{align*}
    Assume that
    \begin{align*}
        \norm{\rho}_{L^\infty(B)} \leq \Lambda, \qquad \norm{\psi}_{L^\infty(\del B)} \leq \Psi_0, \qquad \norm{D^\top \psi}_{L^\infty(\del B)} \leq \Psi_1.
    \end{align*}
    Then there exists a constant $C = C(m, R,\Lambda, \Psi_0, \Psi_1)$ such that
    \begin{align*}
        \sup_{\overline B}w \leq C.
    \end{align*}
\end{thm}

In \cite[Theorem 3.5]{BartnikSimon1982} Bartnik and Simon showed that, if the source $\rho$ is bounded, then interior gradient bounds follow from gradient estimates on the boundary. In light of this, it will suffice to prove the following

\begin{thm}\label{thm:boundary_gradient_estimate}
    Under the assumptions of \cref{thm:full_gradient_estimate}, there exists a constant $C = C(m,R,\Lambda,\Psi_1)$ such that
    \begin{align*}
        \sup_{\del B} w \leq C \sqrt{1 + \Psi_0^2}.
    \end{align*}
\end{thm}

For the Dirichlet problem, Bartnik and Simon obtain a boundary gradient estimate by a barrier argument (see \cite[Corollary 3.4]{BartnikSimon1982}). Here we follow a different approach, more similar to that of \cite{SimonSpruck}, based on a Moser iteration.

We first outline the proof of the gradient estimate, deferring the proofs of the auxiliary estimates to the subsequent subsections. Throughout this Section, $C$ will always denote a constant (whose value might vary from line to line) that depends on $R$, $m$, $\Lambda$ and $\Psi_1$. Further dependencies are indicated by subscripts.

\subsection{Proof of \cref{thm:boundary_gradient_estimate}} The proof of \cref{thm:boundary_gradient_estimate} relies on a Moser iteration scheme on the \emph{rotational tilt function}
\begin{align*}
    w_\top \doteq \left( 1 - \frac{r^2}{R^2}|D^\top u|^2\right)^{- \frac{1}{2}}.
\end{align*}
Here $r=|x|$ and $D^\top u$ denotes the gradient tangential to $\del B_r$ for $r>0$. We set $w_\top(o)$ equal to $1$ at the centre of the ball. The function $w_\top$ extends the one defined in \eqref{eq:boundary_graph_metric} to the whole of $\overline B$ and it is the natural function that arises when differentiating the mean curvature equation solved by $u$ along isometries of the ball $B$. Its main features are
\begin{equation}\label{eq:properties_w_top}
    1 \leq w_\top \leq w, \qquad\quad w_\top \leq \left( 1 - \tfrac{r^2}{R^2} \right)^{-\frac{1}{2}} \quad \text{in ${B}$,} \qquad\quad  w = \sqrt{1 + \psi^2}w_\top \quad \text{on $\del B$}.
\end{equation}
By the last one, since $\psi$ is bounded, in order to achieve a full gradient bound for $u$ on the boundary, it suffices to estimate $w_\top$. By testing the equation for $1-w_\top^{-2}$ with $w_\top^{p+1}$, we obtain the global Caccioppoli inequality:

\begin{prop}\label{prop:caccioppoli}
    For every $p \geq 0$
    \begin{align}\label{eq:caccioppoli}
        \int_{B} |Dw_\top^{p/2}|^2 \leq c_mp^2 \Lambda^2 \int_{B} w_\top^p + c_mp\Psi_1 \int_{\del B} w_\top^{p+1}.
    \end{align}
\end{prop}

In order to improve this estimate to something suitable for a Moser iteration, we need to absorb the boundary integral. More precisely we need an estimate of the form
\begin{align}\label{eq:true_caccioppoli}
    \int_{B} |Dw_\top^{p/2}|^2 \leq Cp^\gamma \int_{B} w_\top^p, \qquad \text{for every $p \geq 2$},
\end{align}
for constants $C$ and $\gamma$ independent of $p$. To this end, we recall the following trace inequality
\begin{align}\label{eq:eps_trace}
    \norm{f}_{L^q(\del B)}^2 \leq \eps \norm{Df}_{L^2(B)}^2 + C_{q}\eps^{-\frac{\theta}{1 - \theta}} \norm{f}_{L^2(B)}^2, \qquad \theta \doteq \tfrac m2 - \tfrac{m-1}{q},
\end{align}
valid for every $\eps \in (0,1]$, $f \in W^{1,2}(B)$ and $2<q<q_*$
where
\begin{align*}
    q_* \doteq 
    \begin{cases}
        \frac{2(m-1)}{m-2} & \qquad \text{if $m \geq 3$,} \\
        \infty & \qquad \text{if $m =2$,}
    \end{cases}
\end{align*}
is the critical trace exponent. This trace inequality follows from Maz'ya's multiplicative inequality \cite[Sec.~1.4.7, Cor.~2]{Mazya2011}, applied to surface measure on $\partial B$, and Young's inequality.
% Inequality \eqref{eq:eps_trace} follows from the calssical Sobolev trace inequality
% \begin{align}
%     \norm{f}_{L^q(\del B)}^2 \leq C_{m,R}\norm{f}_{W^{1,2}(B)}^2, \qquad \forall 1 \leq q \leq q_*
% \end{align}
% Interpolating this latter inequality with $q = q_*$ and the simple estimate
% \begin{align*}
%     \norm{f}_{L^2(\del B)} \leq C_{R,m} \norm{f}_{L^2(B)} \norm{f}_{W^{1,2}(B)}
% \end{align*}
% (see for instance the proof of \cite[Theorem 1.5.1.10]{Grisvard2011}) yields, for every $2 < q < q_*$,
% \begin{align*}
%     \norm{f}_{Lì2(\del B)} \leq C_{R,m,p} \norm{f}_{L^2(B)}^{2(1 - \theta)} \norm{f}_{W^{1,2}(B)}^{2\theta}.
% \end{align*}
% The trace \eqref{eq:eps_trace} follows then by applying Young's inequality.

For any $s > 1$ H\"older's inequality applied to the boundary integral in \eqref{eq:caccioppoli} gives
\begin{align*}
    \int_{\del B} w_\top^{p + 1} \leq \left(\int_{\del B} w_\top^s \right)^{\frac 1s} \left( \int_{\del B} w_{\top}^{\frac{ps}{s-1}} \right)^{\frac{s-1}{s}} = \norm{w_\top}_{L^s(\del B)} \norm{f}_{L^q(\del B)}^2, \qquad q \doteq \frac{2s}{s-1}, \quad f \doteq w_\top^{p/2}.
\end{align*}
The second integral can be estimated by applying \eqref{eq:eps_trace} to $f$ provided that $q < q_*$, that is, $s > m - 1$. In order for this estimate to be effective, we thus need to control $\norm{w_\top}_{L^s(\del B)}$ for some $s > m - 1$. In \cref{subsec:L^s_estimates} we exploit the decaying property \eqref{eq:properties_w_top} of $w_\top$ (cf. \cref{lem:trace_ball}) to achieve much more:

\begin{lem}\label{lem:w_top_Ls_estimate}
    For every $s\geq 2$ there exists a constant $C_s$ such that
    \begin{align}\label{eq:Ls_bounds}
        \int_B |D w_\top^{s/2}|^2 + \int_B w_\top^s + \int_{\del B} w_\top^{s+1} \leq C_s.
    \end{align}
    In particular
    \begin{align}\label{eq:w_top_Ls_estimate}
        \norm{w_\top}_{L^m(\del B)} \leq C_m.
    \end{align}
\end{lem}

Applying the trace inequality \eqref{eq:eps_trace} to $f = w_\top^{p/2}$, using \eqref{eq:w_top_Ls_estimate} and substituting in \eqref{eq:caccioppoli} gives
\begin{align*}
    \int_B |D w_\top^{p/2}|^2 \leq Cp^2 \int_B w_\top^p + \eps C_mCp \int_B |Dw_\top^{p/2}|^2 + \eps^{-\frac{\theta}{1 - \theta}}C_m Cp  \int_B w_\top^p.
\end{align*}
We set $\eps = \frac{1}{1 + 2C_mCp}$, so that $\eps C_m Cp \leq \frac12$ and the energy can be absorbed in the left-hand side:
\begin{align*}
    \int_B |D w_\top^{p/2}|^2 \leq Cp \left(p + (1 + 2C_m Cp)^{\frac{\theta}{1 - \theta}}C_m C \right) \int_B w_\top^p.
\end{align*}
This gives \eqref{eq:true_caccioppoli} with $\gamma = \max\set{2,1 + \tfrac{\theta}{1 - \theta}}$.

Combining the $L^2$-Sobolev inequality and \eqref{eq:true_caccioppoli} we have, for every $p\geq2$,
\begin{align*}
    \left(\int_B w_\top^{\chi p} \right)^{\frac{1}{\chi}} \leq C \left( \int_B w_\top^p + \int_B |Dw_\top^{p/2}|^2 \right) \leq Cp^\gamma \int_B w_\top^p, \qquad \chi \doteq \frac{m}{m-2} \bm 1_{\set{m\geq 3}} + 2 \, \bm 1_{\set{m=2}}.
\end{align*}
Taking the $p$-th root and iterating $N$ times starting from $p = 2$, letting $p_j \doteq 2\chi^j$ we obtain
\begin{align}\label{eq:moser_iteration}
    \left(\int_B w_\top^{p_N} \right)^{\frac{1}{p_N}} \leq \prod_{j=0}^{N-1} \left(C p_j^\gamma \right)^{\frac{1}{p_j}} \left( \int_B w_\top^2 \right)^\frac12.
\end{align}
The product converges because the series
\begin{align*}
    \log \left(\prod_{j=0}^{\infty}\left(C p_j^\gamma \right)^{1/p_j} \right) = \sum_{j=0}^{\infty} \frac{\log C + \gamma \log 2 + j\gamma \log \chi } {2\chi^j}
\end{align*}
is finite. As a consequence taking the limit $N \to \infty$ in \eqref{eq:moser_iteration} gives
\begin{align*}
    \norm{w_\top}_{L^\infty(B)} \leq C \norm{w_\top}_{L^2(B)}.
\end{align*}
\cref{lem:w_top_Ls_estimate} also gives the $L^2$ estimate
\begin{align}\label{eq:w_top_L2_estimate}
    \norm{w_\top}_{L^2(B)} \leq C_2,
\end{align}
hence $w_\top$ is bounded by a constant that only depends on $R$, $m$, $\Lambda$ and $\Psi_1$. Since $u \in C^2(\overline{B})$, this estimate holds on the boundary as well and by \eqref{eq:properties_w_top} the estimate follows.

% \begin{rmk}
%     Without loss of generality we can assume $R = 1$. Indeed, the rescaled function $\tilde u (x) \doteq \frac1R u(Rx)$ defined on $B_1$ solves \eqref{eq:problem_Rm} with data
%     \begin{align*}
%         \tilde\rho (x) \doteq R \, \rho(Rx) \qquad \tilde\psi(x) = \psi(Rx)
%     \end{align*}
%     that satisfy $\norm{\tilde\rho}_{L^\infty(B_1)} = R \norm{\rho}_{L^\infty(B_R)}$, $\norm{\tilde \psi}_{L^\infty(\del B_1)} = \norm{ \psi}_{L^\infty(\del B_R)}$ and $\norm{D^\top\tilde \psi}_{L^\infty(\del B_1)} = R\norm{D^\top \psi}_{L^\infty(\del B_R)}$.
% \end{rmk}

\subsection{The Caccioppoli estimate}

In this paragraph we will prove the Caccioppoli estimate in \cref{prop:caccioppoli}. As anticipated, this is obtained by differentiating the mean curvature equation along rotations of the ball $B$ and testing the resultant PDE with $w_\top^{p+1}$.

For every $i < j$ let $Q^\eps_{ij} \in SO(m)$ be the counter-clockwise rotation of angle $\eps>0$ in the plane $(i,j)$ (so that $Q^\eps_{ij}\del_k = \del_k$ for every $k \neq i,j$) and let
\begin{align*}
    u^\eps_{ij}(x) \doteq u(Q^\eps_{ij}x), \qquad z_{ij}(x) \doteq \frac{d}{d\eps}\bigg\rvert_{\eps = 0} u(Q^\eps_{ij}(x)).
\end{align*}
Alternatively, letting
\begin{align*}
    X_{ij} \doteq x_i \del_j - x_j \del_i
\end{align*}
be the infinitesimal generator of $Q^\eps_{ij}$, we have $z_{ij} = X_{ij}(u)$. The rotational derivatives $z_{ij}$ solve the linearized equation
\begin{align}\label{eq:rotation_equation}
    - \div(wg^{-1}Dz_{ij}) = X_{ij}(\rho) \quad \text{in $B_R$}, \qquad -wg^{-1}Dz_{ij} \cdot n = X_{ij}(\psi) \quad \text{on $\del B_R$}
\end{align}
where
\begin{align*}
    g^{-1} \doteq \id + w^2 Du \otimes Du.
\end{align*}
This is the inverse graph metric in Cartesian coordinates, see \cref{sec:preliminaries}. Since the eigenvalues of $g^{-1}$ are $w^2$ and $1$ (with multiplicity $1$ and $m - 1$ respectively), we have
\begin{align}\label{eq:B_bound}
    g^{-1} \geq \id.
\end{align}
Consider now
\begin{align*}
    z^2 \doteq \frac{1}{R^2}\sum_{i < j} z_{ij}^2 = \frac{r^2}{R^2} |D^\top u|^2, \qquad D^\top u \doteq Du - \left(Du \cdot \frac{x}{r} \right) \frac{x}{r}, \qquad r \doteq |x|.
\end{align*}
Here $D^\top u$ is defined for $r>0$, while $z^2$ is defined at $x=0$ by the sum above. By \eqref{eq:rotation_equation}, the function $z^2$ solves
\begin{equation}\label{eq:PDE_z}
    \begin{cases}
        \div(wg^{-1}Dz^2) = \dfrac{2}{R^2}\sum_{i<j}wg^{-1}(Dz_{ij},Dz_{ij}) - \dfrac{2}{R^2}\sum_{i<j}z_{ij}X_{ij}(\rho) & \qquad \text{in $B_R$},\\
        -wg^{-1}Dz^2\cdot n = 2D^\top u \cdot D^\top\psi & \qquad \text{on $\partial B_R$}.
    \end{cases}
  %   \begin{aligned}
  %   \div(wg^{-1}Dz)
  %   &= \frac{2}{R^2}\sum_{i<j}wg^{-1}(Dz_{ij},Dz_{ij}) - \frac{2}{R^2}\sum_{i<j}z_{ij}X_{ij}(\rho)
  %   &&\text{in }B_R,\\
  % wg^{-1}Dz\cdot n
  % &=2D^\top u \cdot D^\top\psi
  %   &&\text{on }\partial B_R.
  %   \end{aligned}
\end{equation}
% Consider the \emph{rotational tilt function}
% \begin{align*}
%     w_\top \doteq \frac{1}{\sqrt{1 - \tfrac{r^2}{R^2} |D^\top u|^2}} = \frac{1}{\sqrt{1 - z}}.
% \end{align*}

\begin{proof}[Proof of \cref{prop:caccioppoli}]
    The case $p=0$ is immediate, so assume $p>0$.
    % Observe that
    % \begin{equation}\label{eq:properties_w_top}
    %     1 \leq w_\top \leq w \quad \text{in $\overline{B_R}$}, \qquad w_\top \leq \left( 1 - \tfrac{r^2}{R^2} \right)^{-\frac{1}{2}} \quad \text{in $B_R$},
    % \end{equation}
    % while
    % \begin{align}\label{eq:w_top_at_boundary}
    %     w = \sqrt{1 + \psi^2}w_\top \qquad \text{on $\del B_R$}.
    % \end{align}
    Multiplying \eqref{eq:PDE_z} by $\varphi/2$, where $\varphi \in C^1(\overline{B_R})$ is nonnegative, and integrating by parts we get
    \begin{gather*}
        \frac{1}{R^2}\sum_{i < j} \int_{B_R} \varphi wg^{-1}(Dz_{ij},Dz_{ij}) + \frac12\int_{B_R} wg^{-1}(Dz^2,D\varphi) = -\int_{\del B_R} \varphi D^\top u \cdot D^\top \psi + \frac{1}{R^2}\sum_{i<j} \int_{B_R} \varphi z_{ij}X_{ij}(\rho).
    \end{gather*}
    Since the fields $X_{ij}$ are tangential and divergenceless we have
    \begin{align*}
        \sum_{i<j} \int_{B_R} \varphi z_{ij} X_{ij}(\rho) = - \sum_{i<j} \int_{B_R} \rho \div(\varphi z_{ij}X_{ij}) = - \sum_{i<j} \int_{B_R} \rho z_{ij} X_{ij}(\varphi) - \sum_{i<j} \int_{B_R} \rho \varphi X_{ij}(z_{ij})
    \end{align*}
    obtaining
    \begin{gather*}
        \frac{1}{R^2}\sum_{i < j} \int_{B_R} \varphi wg^{-1}(Dz_{ij},Dz_{ij}) + \frac12\int_{B_R} wg^{-1}(Dz^2,D\varphi) \leq \qquad\qquad\qquad \\ \qquad\qquad\qquad \leq \Psi_1\int_{\del B_R} \varphi - \frac{1}{R^2}\int_{B_R} \rho r^2 D^\top u \cdot D\varphi - \frac{1}{R^2}\sum_{i<j} \int_{B_R} \rho \varphi X_{ij}(z_{ij}),
    \end{gather*}
    where we also used that
    \begin{align*}
        \sum_{i<j} z_{ij}X_{ij} = r^2 D^\top u, \qquad |D^\top u \cdot D^\top \psi| \leq \Psi_1.
    \end{align*}
    We can estimate $X_{ij}(z_{ij})$ right away by Cauchy--Schwarz, \eqref{eq:B_bound} and $\sum_{i<j}|X_{ij}|^2 = (m-1)r^2$:
    \begin{align*}
        \Big|\sum_{i<j} X_{ij}(z_{ij}) \Big|^2 \leq \Bigg(\sum_{i<j} |X_{ij}|^2 \Bigg) \Bigg(\sum_{i<j} |Dz_{ij}|^2 \Bigg) \leq r^2(m - 1) \ \frac{1}{w}\sum_{i<j} wg^{-1}(Dz_{ij}, Dz_{ij}).
    \end{align*}
    Hence, after applying Young's inequality
    \begin{gather*}
        \frac{1}{2R^2}\sum_{i < j} \int_{B_R} \varphi wg^{-1}(Dz_{ij},Dz_{ij}) + \frac12\int_{B_R} wg^{-1}(Dz^2,D\varphi) \leq \qquad\qquad\qquad \\ \qquad\qquad\qquad \leq \Psi_1\int_{\del B_R} \varphi - \frac{1}{R^2}\int_{B_R} \rho r^2 D^\top u \cdot D\varphi + \frac{m - 1}{2R^2}\Lambda^2\int_{B_R} \frac{\varphi r^2}{w} .
    \end{gather*}
    Choosing the test function
    \begin{align*}
        \text{$\varphi = w_\top^{p+1}$,} \qquad \text{$D\varphi = \frac{p + 1}{2} w_\top^{p+3} Dz^2$}
    \end{align*}
    and using the positivity of $wg^{-1}$ to discard the term with $wg^{-1}(Dz_{ij}, Dz_{ij})$ and estimating $r^2/R^2 \leq 1$, $w_\top \leq w$ we arrive at
    \begin{align*}
        \frac{p+1}{4}\int_{B_R} w_\top^{p+3} wg^{-1}(Dz^2,Dz^2) \leq \Psi_1 \int_{\del B_R} w_\top^{p+1} + \frac{p+1}{2}\int_{B_R} w_\top^{p+3}|\rho|\, |D^\top u \cdot Dz^2| + \tfrac{m - 1}{2}\Lambda^2\int_{B_R} w_\top^p.
    \end{align*}
    In order to absorb the term $|D^\top u \cdot Dz^2|$, we apply the Cauchy--Schwarz inequality with respect to $wg^{-1}$:
    \begin{align*}
        |D^\top u \cdot Dz^2|  \leq \sqrt{w^{-1}g(D^\top u, D^\top u)} \, \sqrt{wg^{-1}(Dz^2,Dz^2)}.
    \end{align*}
    Since $w^{-1}g = w^{-1} \left(\id - \, du \otimes du \right)$ we have
    \begin{align*}
        w^{-1}g(D^\top u, D^\top u) = \frac{|D^\top u|^2 (1 - |D^\top u|^2)}{w} \leq \frac{1}{ww_\top^2},
    \end{align*}
    hence by applying Young's inequality once again we have
    \begin{align*}
        \tfrac{p+1}{8}\int_{B_R} w_\top^{p+3} wg^{-1}(Dz^2,Dz^2) \leq \Psi_1 \int_{\del B_R} w_\top^{p+1} + \tfrac{m + p}{2}\Lambda^2\int_{B_R} w_\top^p.
    \end{align*}
    Finally we use $Dz^2= 2w_\top^{-3} Dw_\top$ and \eqref{eq:B_bound} to estimate
    \begin{align*}
        \frac{p+1}{8}w_\top^{p + 3} wg^{-1}(Dz^2,Dz^2) \geq \frac{p+1}{2} w_\top^{p - 2}|Dw_\top|^2 = \frac{p+1}{2} \left| w_\top^{\frac{p}{2} - 1}Dw_\top \right|^2 = \frac{2(p+1)}{p^2} \left| D(w_\top^{\frac p2}) \right|^2
    \end{align*}
    and the estimate follows using $m \geq 2$ and $p>0$.
\end{proof}

\subsection{$L^s$-estimates}\label{subsec:L^s_estimates}
This Subsection is devoted to proving \cref{lem:w_top_Ls_estimate}. It follows from the Caccioppoli estimate \cref{prop:caccioppoli} and the following two lemmas. The first one relies on the decaying property
% $w_\top \leq \left(1 - \tfrac{r^2}{R^2} \right)^{-1/2}$
(the second in \eqref{eq:properties_w_top}) which implies that if $w_\top$ takes high values on the boundary, then it must decrease rapidly towards the interior, entailing a cost in terms of the Dirichlet energy. We stress that the following Lemmas do not require either $w_\top$ or $u$ to solve a PDE.

\begin{lem}\label{lem:trace_ball}
    Under the assumptions above, for every $s \geq 2$ and $\eps > 0$ it holds
    \begin{align}\label{eq:boundary_L^p+1}
        \int_{\del B_R} w_\top^{s+1} &\leq C_{s,\eps} + \eps \int_{B_R} |Dw_\top^{\frac s2}|^2.
    \end{align}
\end{lem}
\begin{proof}
    For a chosen $\omega \in \S^{m-1}$ consider the function $h(r) \doteq w_\top^{\frac s2}(r\omega)$. Letting $\rho \doteq R - r$ be the distance from the boundary,
    % since $1 - \frac{r^2}{R^2} = \frac{(R + r)(R - r)}{R^2} \geq \frac{\rho}{R}$,
    the decaying property in \eqref{eq:properties_w_top} implies
    \begin{align*}
        h(r)^{\frac 2s} \leq \sqrt{\frac{R}{\rho}},
    \end{align*}
    hence, at $\rho_* = R - R_* \doteq \frac{4R}{h(R)^{4/s}}$ the function $h^{\frac2s}$ is at most half its boundary value, provided $h(R)^{2/s} \geq 4$ (which ensures $\rho_* \leq R/4$). In this case
    \begin{align*}
        h(R) - h(R_*) \geq (1 - 2^{-s/2})h(R) \geq \tfrac{1}{2} h(R).
    \end{align*}
    Now, by the fundamental Theorem of calculus and the Cauchy--Schwarz inequality we have
    \begin{align*}
        \tfrac{1}{2}h(R)
        &\leq |h(R) - h(R_*)|
        \leq \int_{R_*}^R |h'(r)| \, dr \\
        &\leq (R - R_*)^{1/2} \left(\int_{R_*}^R |h'|^2 \, dr \right)^{1/2}
        = \frac{2\sqrt R}{h(R)^{2/s}} \left(\int_{R_*}^R |h'|^2 \, dr \right)^{1/2},
    \end{align*}
    that is, estimating the radial derivative with the full gradient of $w_\top^{s/2}$,
    \begin{align*}
        w_\top^{s + 2}(R\omega) \leq 16 R \int_{R_*}^R |Dw_\top^{\frac s2}(r\omega)|^2 \, dr, \qquad \forall \omega \in \S^{m-1} \ : \ w_\top(R\omega) \geq 4.
    \end{align*}
    Now we integrate on $E \doteq \set{\omega \in \S^{m-1} \ | \ w_\top(R\omega) \geq 4}$. Since
    \begin{align*}
        R_* = R \left( 1 - \tfrac{4}{w_\top(R\omega)^2} \right) \geq \tfrac{R}{2} \qquad \forall \omega \in E,
    \end{align*}
    we have
    \begin{align*}
        R^{m - 1}\int_E w_\top^{s+2}(R\omega) \, d\omega
        &\leq 16 R^m \int_E \int_{R_*}^R |Dw_\top^{\frac s2}(r\omega)|^2 \, dr \, d\omega \\
        &\leq 2^{m+3} R \int_E \int_{R_*}^R |Dw_\top^{\frac s2}(r\omega)|^2 r^{m - 1} \, dr \, d\omega.
    \end{align*}
    The integral on $E^c$ is estimated easily:
    \begin{align*}
        R^{m-1} \int_{E^c} w_\top^{s+2}(R\omega) \, d\omega \leq 4^{s+2} |\del B_R|,
    \end{align*}
    and summing up the two contributions gives
    \begin{align*}
        \int_{\del B_R} w_\top^{s+2} \leq 2^{m+3} R \int_{B_R} |D(w_\top^{s/2})|^2 + 4^{s+2}|\del B_R|.
    \end{align*}
    The claim follows by applying the inequality $t^{s+1} \leq c t^{s+2} + C_{s,c}$ for $t \geq 0$, with $c = \frac{\eps}{2^{m+3}R}$.
\end{proof}

\begin{lem}
    Under the assumptions above, for every $s \geq 2$ and $\eps > 0$ it holds
    \begin{align}\label{eq:interior_Lp}
        \int_{B_R} w_\top^s &\leq C_{s,\eps} + \eps \int_{B_R} |Dw_\top^{\frac s2}|^2.
    \end{align}
\end{lem}
\begin{proof}
    Let $\delta \in (0,\tfrac 12)$. In the interior region $B_{(1 - \delta)R}$, by \eqref{eq:properties_w_top} we have $w_\top \leq \delta^{-1/2}$, hence
    \begin{align*}
        \int_{B_{(1 - \delta)R}} w_\top^s \leq \delta^{-\frac s2}|B_R|.
    \end{align*}
    To treat the complement, fix $\omega \in \S^{m-1}$. By the fundamental Theorem of calculus applied to $w_\top^{\frac s2}(\omega t)$ and the Cauchy--Schwarz inequality, for $(1-\delta)R \leq r \leq R$ we have
    \begin{align*}
        w_\top^s(\omega r)
        &\leq 2 w_\top^s(\omega (1 - \delta)R) + 2\delta R \int_{(1 - \delta)R}^R |\partial_t w_\top^{\frac s2}(\omega t)|^2 \, dt \\
        &\leq 2\delta^{-\frac{s}{2}} + 2\delta R \int_{(1 - \delta)R}^R |\partial_t w_\top^{\frac s2}(\omega t)|^2 \, dt.
    \end{align*}
    Hence, integrating in polar coordinates in the region $B_R\setminus B_{(1-\delta)R}$ and using $t \geq (1-\delta)R \geq R/2$, we have
    \begin{align*}
        \int_{B_R\setminus B_{(1-\delta)R}} w_\top^s
        &\leq 2\delta^{-\frac{s}{2}}|B_R| + \frac{2\delta^2 R^2}{(1-\delta)^{m-1}} \int_{B_R} |Dw_\top^{\frac s2}|^2 \\
        &\leq 2\delta^{-\frac{s}{2}}|B_R| + 2^m\delta^2 R^2 \int_{B_R} |Dw_\top^{\frac s2}|^2.
    \end{align*}
    Summing the two contributions and choosing $\delta$ small enough that $2^m\delta^2 R^2 \leq \eps$ yields \eqref{eq:interior_Lp}.
\end{proof}

Applying the last two estimates in \eqref{eq:caccioppoli} with $p=s$ and choosing $\eps>0$ small enough to absorb the gradient terms gives \eqref{eq:Ls_bounds}.

\section{Higher integrability of the tilt function}

In this Section, we will establish integral bounds for the tilt function of smooth, strictly space-like classical solutions to \eqref{eq:problem_Rm} on a ball in terms of the data $\rho$ and $\psi$. 

\subsection{Second fundamental form estimates}

Let $F:\overline\Omega \to M$ be a $C^3$-regular, strictly space-like graph as in \cref{sec:preliminaries} and let $H = -\rho$ be its mean curvature. In the present Subsection, we establish a priori global integral bounds on the second fundamental form and on the square of the tilt function. These estimates ultimately rely on a couple of geometric facts: the Codazzi equation (see e.g. \cite[Section 9.2]{AliasMastroliaRigoli})
\begin{align*}
    \div A = - \nabla \rho
\end{align*}
and the fact that $\del_0$ is a Killing vector field of $\R \times \Omega\subset \Mink$.

Call $X^M\in \X(\overline\Omega)$ the tangential component of a vector field $X \in \X(\Mink)$, defined by $X = F_*X^M - \langle X,N\rangle N.$ In particular, $\del_0$ decomposes as
\begin{align*}
    \del_0 = F_*\del_0^M + wN, \qquad \del_0^M = - \nabla u
\end{align*}
where $u \doteq - \langle F,\del_0\rangle$ is the height function of the graph, as usual. Since $\del_0$ is parallel, differentiating $w = -\langle N,\del_0\rangle$ gives
\begin{align*}
    \nabla w = A\nabla u.
\end{align*}
The Codazzi equation can be rephrased as
\begin{align*}
    \div_g T = 0, \qquad T \doteq A + \rho\id.
\end{align*}
Notice that $T$ is minus the first Newton tensor of $A$. Thus
\begin{align}\label{eq:divergence_X}
    \div_g X = w \left( |\II|^2 - \rho^2 \right), \qquad X \doteq -T(\del_0^M) = \nabla w + \rho\nabla u.
\end{align}
This vector field is geometrically natural, and it is designed to feature second derivatives of $u$ without involving third ones. A closely related vector field was used in \cite{BIMM} to obtain similar estimates.

\begin{lem}\label{lem:second_identity}
    Let $u \in C^3(\overline B)$ be a zero-mean function such that $|Du| < 1$ on $\overline B$ and let
    \begin{align*}
        - \div(wDu) = \rho \quad \text{in $B$}, \qquad - w Du \cdot n = \psi \quad \text{on $\del B$}.
    \end{align*}
    Then
    \begin{align}\label{eq:second_identity_split}
        \int_B \left(|\II|^2 - \rho^2\right) \, dx + \frac1R \int_{\del B}\left\{ w^2|D^\top u|^2 + (m-1)\psi^2 \right\} \, d\Haus^{m-1} + 2 \int_{\del B} wD^\top u \cdot D^\top \psi \, d\Haus^{m-1} = 0. 
    \end{align}
\end{lem}
\begin{proof}
    Integrating by parts \eqref{eq:divergence_X} in the graph metric and recalling \eqref{eq:volumes}, we obtain
    \begin{align*}
        \int_B |\II|^2 \, dx = \int_B \rho^2 \, dx + \int_{\del B} g(X,\nu)\, d\sigma.
    \end{align*}
    Here $d\sigma$ is the Riemannian volume of the boundary of the graph, as in \cref{sec:preliminaries}. We rewrite the boundary integral in terms of the Euclidean geometry. %First of all notice that
    % \begin{align*}
    %     X = Dw - w^3 \Delta u \, Du.
    % \end{align*}
    Combining \eqref{eq:conormal} and \eqref{eq:boundary_volumes} we obtain
    \begin{align}\label{eq:flux_formula}
        g(\nabla f, \nu) \, d\sigma = \left( \frac{f_n}{w} - \psi Du \cdot Df \right) \, d\Haus^{m-1}, \qquad f_n \doteq \del_n f,
    \end{align}
    for every $f \in C^1(\overline B)$. In particular
    \begin{align*}
        g(\nabla w, \nu) \, d\sigma = \left( \frac{w_n}{w} - \psi \, Du \cdot Dw \right) \, d\Haus^{m-1}, \qquad g(\nabla u, \nu) \, d\sigma = \left(\frac{u_n}{w} - \psi \, |Du|^2\right) \, d\Haus^{m-1} = - \psi \, d\Haus^{m-1}.
    \end{align*}
    Hence, recalling also that $-\rho = w\Delta u + Dw \cdot Du$ and $Dw = w^3D^2u(Du)$, we have
    \begin{align*}
        g(X,\nu) \, d\sigma = \left(D^2u(Du,n) - u_n \Delta u \right) w^2 \, d\Haus^{m-1} = \left\{\del_n \ln w + w \psi \Delta u \right\} \, d\Haus^{m-1}.
    \end{align*}
    Next we notice that on the sphere $\del B$ one has
    \begin{align*}
        D^2u(D^\top u, n) = D^\top u \cdot D^\top u_n - \frac{|D^\top u|^2}{R}, \qquad \Delta u = \Delta_{\del B} u +  D^2u(n,n) + \frac{m-1}{R} u_n,
    \end{align*}
    so
    \begin{align*}
        \int_{\del B} g(X,\nu) \, d\sigma &= \int_{\del B} \left\{ w^2D^\top u \cdot D^\top(u_n) - u_n w^2 \Delta_{\del B} u - \frac{w^2|D^\top u|^2}{R} - \frac{m-1}{R}\psi^2 \right\} \, d\Haus^{m-1} \\
        &= \int_{\del B} \left\{ -wD^\top u \cdot D^\top \psi + \psi \div_{\del B} (wD^\top u) - \frac{w^2|D^\top u|^2}{R} - \frac{m-1}{R}\psi^2 \right\} \, d\Haus^{m-1}.
    \end{align*}
    Integrating by parts we have
    \begin{align*}
        \int_{\del B} g(X,\nu) \, d\sigma &= \int_{\del B} \left\{ - 2wD^\top u \cdot D^\top \psi - \frac{w^2|D^\top u|^2}{R} - \frac{m-1}{R}\psi^2 \right\} \, d\Haus^{m-1}.
    \end{align*}
    This proves the identity.
\end{proof}

The integral identity in \cref{lem:second_identity} allows us to obtain integral bounds for the tilt function and the second derivatives of $u$ in terms of the $L^2$ norm of the mean curvature, and of the $H^1$ norm of the capillary datum.

\begin{lem}\label{lem:L^2_w_estimates}
    Under the assumptions of \cref{lem:second_identity} we have
    \begin{align}
        \norm{wD^2u}_{L^2(B)}^2 + \norm{D \ln w}_{L^2(B)}^2 + \norm{Dw \cdot Du}_{L^2(B)}^2 \leq \norm{\rho}_{L^2(B)}^2 + R \norm{D^\top \psi}_{L^2(\del B)}^2
    \end{align}
    and
    \begin{align}
        \norm{w}_{L^2(B)}^2 \leq (m+2)|B_R| + 4R^2 \left( \norm{\rho}_{L^2(B)}^2 + 2R \norm{D^\top \psi}_{L^2(\del B)}^2 \right).
    \end{align}
\end{lem}
\begin{proof}
    The first estimate follows immediately from \eqref{eq:second_identity_split} and
    \begin{align*}
        |\II|^2 = w^2 |D^2u|^2 + 2w^4 |D^2u(Du)|^2 + w^6 D^2u(Du,Du)^2 = w^2|D^2u|^2 + 2 |D\ln w|^2 + (Dw \cdot Du)^2.
    \end{align*}

    For the second one, we apply Young's inequality to the last integral in \eqref{eq:second_identity_split} 
    % \begin{align*}
    %     \int_B |\II|^2 + \frac{1}{2R} \int_{\del B} w^2|D^\top u|^2 + \frac{m-1}{R} \int_{\del B} \psi^2 \leq \int_B \rho^2 + 2R \int_{\del B} |D^\top \psi|^2 
    % \end{align*}
    together with the identity 
    \begin{align*}
        w^2|D^\top u|^2 = w^2|Du|^2 - w^2 u_n^2 = w^2 - \psi^2 - 1
    \end{align*}
    and $m \geq 2$ to obtain
    \begin{align}\label{eq:second_intermediate_estimate}
        \int_B |\II|^2 + \frac{1}{2R} \int_{\del B} \left[w^2 + (2m-3)\psi^2 \right] \leq \int_B \rho^2 + 2R \int_{\del B} |D^\top \psi|^2 + \frac{|\del B_R|}{2R}.
    \end{align}
    Up to subtracting from $u$ its value in the centre of the ball, we can assume $|u| \leq R$.
    Now, testing the mean curvature equation with $uw$ gives
    \begin{align*}
        \frac 1R\int_B (w^2 - 1) &= - \frac 1R \int_B uw Du\cdot Dw + \frac 1R\int_B \rho \,  uw - \frac 1R\int_{\del B} \psi \, uw \\ 
        & \leq \int_B w |Du \cdot Dw| + \int_B |\rho| \, w +  \int_{\del B} |\psi| \, w \\
        &\leq \frac{1}{2R} \int_B w^2 + R \int_B (Du \cdot Dw)^2 + R \int_B \rho^2 + \frac{1}{2} \int_{\del B} (w^2 + \psi^2)
    \end{align*}
    hence using also $|\II|^2 \geq (Du\cdot Dw)^2$ and \eqref{eq:second_intermediate_estimate}
    \begin{align*}
        \frac{1}{2R^2} \int_B w^2 &\leq \int_B |\II|^2 + \frac{1}{2R} \int_{\del B} (w^2 + \psi^2) + \int_B \rho^2 + \frac{|B_R|}{R^2} \\
        &\leq 2 \int_B \rho^2 + 2R\int_{\del B} |D^\top \psi|^2 + \frac{|B_R|}{R^2} + \frac{|\del B|}{2R}
    \end{align*}
    which proves the claim.
\end{proof}

\subsection{Uniform integrability of the tilt function}

\begin{thm}\label{thm:boundary_energy_estimate}
    Let $u$ be as in \cref{lem:second_identity} and, for every $\delta>0$,
    let $U_\delta \doteq \set{R-\delta<|x|<R}$. Fix
    $\delta_0\in(0,R/2)$. Then, for every $q\in(1,\infty]$, there exists
    a constant
    \begin{align*}
        C=C\left(m,R,q,\norm{\rho}_{L^1(B)},\norm{\psi}_{L^q(\del B)}\right)
    \end{align*}
    such that, for every $0<\delta\leq\delta_0/2$,
    \begin{align*}
        \int_{U_\delta} w\,dx \leq C\left(\delta^{\alpha_q}+\int_{U_{2\delta}}|\rho|\,dx\right),
    \end{align*}
    where
    \begin{align*}
        \alpha_q = \frac{q-1}{4(2q-1)} \bm 1_{\set{1<q<2}} + \frac{1}{12} \bm 1_{\set{2\leq q \leq \infty}}.
    \end{align*}
\end{thm}

% Combining \eqref{thm:boundary_energy_estimate}, the no-light-segment \cref{thm:no_light_segment} and the higher integrability results in \cite{BIMM} we obtain

% \begin{cor}\label{eq:uniform_integrability}
%     Assume $\rho \in $
% \end{cor}

We give the proof of \cref{thm:boundary_energy_estimate} in five steps. 

\subsubsection*{\textbf{Step 1. Reduction to $R = 1$}} \textit{It suffices to prove the Theorem with $R = 1$.} Indeed, if $u$ solves \eqref{eq:problem_Rm} in $B_R$ with data $\rho$ and $\psi$, then the rescaling
\begin{align*}
    \tilde u(x) \doteq \tfrac1R u(Rx), \qquad x \in B_1
\end{align*}
solves \eqref{eq:problem_Rm} in $B_1$ with data
\begin{align*}
    \tilde\rho (x) \doteq R \rho(Rx), \qquad \tilde\psi(\omega) \doteq \psi(R\omega).
\end{align*}
Moreover
\begin{gather*}
    \norm{\tilde\rho}_{L^1(B_1)} = R^{1 - m} \norm{\rho}_{L^1(B_R)}, \qquad \norm{\tilde\psi}_{L^q(\del B_1)} = R^{-\frac{m-1}{q}} \norm{\psi}_{L^q(\del B_R)}, \\
    \int_{U_{2\delta/R}(B_1)}|\tilde\rho| \, dx = R^{1-m} \int_{U_{2\delta}(B_R)}|\rho| \, dx.
\end{gather*}

\subsubsection*{\textbf{Step 2. Basic local estimates}}
\textit{Fix $\delta_0\in(0,1/2)$. Here we show that there exists a constant $C=C(m,q,\norm{\rho}_{L^1(B)}, \norm{\psi}_{L^q(\del B)})$ such that, for every $0<\delta\leq\delta_0/2$,}
\begin{align}\label{eq:easy_energy_integral_bounds}
    \int_B w^2 \, d\mu \leq C, \qquad
    \int_{U_\delta} |\nabla r|^2 \, d\mu \leq C\delta, \qquad
    \int_{U_\delta} w|\nabla r| \, d\mu \leq C\sqrt{\delta},
\end{align}
\textit{and that, for every $\eta\in C^\infty(\del B)$ with $0\leq \eta\leq 1$, extended radially to the interior, we have}
\begin{align}\label{eq:ikari}
    \int_{U_\delta} \eta w^2 \, d\mu
    \leq C\delta+ C\norm{D^\top\eta}_{L^\infty(\del B)}
    \sqrt{\delta}+\int_{\del B}\eta\cosh\gamma\,d\sigma + \int_{U_{2\delta}} |\rho| \, dx.
\end{align}

The first estimate in \eqref{eq:easy_energy_integral_bounds} follows by testing \eqref{eq:mean_curv_intrinsic} with $u$ itself and using \eqref{eq:intrinsic_normal_derivatives} together with the first identity in \eqref{eq:useful_identities}:
\begin{align*}
    \int_B(w^2-1)\,d\mu
    =\int_B|\nabla u|^2\,d\mu
    =\int_B u\rho\,dx-\int_{\del B}u\psi\,d\omega.
\end{align*}
Indeed, by compatibility we may subtract $u(0)$ from $u$ on the right-hand side. Since $|u-u(0)|\leq 1$ and $d\mu=w^{-1}\,dx\leq dx$, the estimate follows.

The second estimate is obtained by testing the $g$-Laplacian of the Euclidean distance $r$ from the centre of the ball:
\begin{align}\label{eq:laplacian_r}
    \Delta_g r
    =\frac{w^2-|\nabla r|^2}{r}+\frac{m-1}{r}
    -\frac{\rho}{w}\langle\nabla u,\nabla r\rangle.
\end{align}
This can be seen by taking the Laplacian in $r^2=\sum_i x_i^2$, tracing $g=\delta_{ij}\,dx^i\otimes dx^j-du^2$ with respect to $g$ and using the third identity in \eqref{eq:useful_identities} together with the last one in \eqref{eq:mean_curv_intrinsic}. As a test function we choose
\begin{align*}
    \chi_\delta(r)\doteq\left(1-\tfrac{1-r}{\delta}\right)_+
\end{align*}
and integrate by parts the divergence of $\chi_\delta\nabla r$. By \eqref{eq:intrinsic_normal_derivatives} and \eqref{eq:laplacian_r} we obtain
\begin{align*}
    &\frac{1}{\delta}\int_{U_\delta}|\nabla r|^2\,d\mu
    +\int_{U_\delta}\frac{\chi_\delta}{r}
        \left(w^2-|\nabla r|^2+m-1\right)\,d\mu=\\
    &\qquad=\int_{\del B}\cosh\gamma\,d\sigma
    +\int_{U_\delta}\chi_\delta\frac{\rho}{w}
        \langle\nabla u,\nabla r\rangle\,d\mu.
\end{align*}
Since $|\nabla r|\leq w$, the second integral on the left-hand side is non-negative and can be discarded. Moreover, since also $|\nabla u|\leq w$,
\begin{align}\label{eva_00}
    \left|\frac{\rho}{w}\langle\nabla u,\nabla r\rangle\right|
    \leq |\rho|w.
\end{align}
As $w\,d\mu=dx$, the last integral is bounded in absolute value by $\norm{\rho}_{L^1(B)}$. Since $\cosh\gamma=\sqrt{1+\psi^2}$ and $d\sigma\leq d\omega$, this gives the second estimate in \eqref{eq:easy_energy_integral_bounds}. The third follows from the first two by H\"older's inequality. These three estimates hold, in fact, for every $0<\delta\leq\delta_0$.

Estimate \eqref{eq:ikari} is obtained by integrating by parts the vector field $\eta\xi_\delta\nabla r$, where $\xi_\delta$ is a smooth, non-decreasing radial cut-off such that
\begin{align*}
    0\leq\xi_\delta\leq 1,\qquad
    \xi_\delta=1 \quad\text{in $U_\delta$},\qquad
    \xi_\delta=0 \quad\text{in $B\setminus U_{2\delta}$}.
\end{align*}
We obtain
\begin{align*}
    &\int_{U_{2\delta}}\frac{\eta\xi_\delta}{r}
        \left(w^2-|\nabla r|^2+m-1\right)\,d\mu
    +\int_B\eta\xi_\delta'|\nabla r|^2\,d\mu=\\
    &\qquad=\int_{\del B}\eta\cosh\gamma\,d\sigma
    +\int_B\eta\xi_\delta\frac{\rho}{w}
        \langle\nabla u,\nabla r\rangle\,d\mu
    -\int_B\xi_\delta\langle\nabla r,\nabla\eta\rangle\,d\mu.
\end{align*}
Both integrands on the left-hand side are non-negative. Restricting the first integral to $U_\delta$, where $\xi_\delta=1$ and $1/r\geq1$, we see that the left-hand side is bounded below by
\begin{align*}
    \int_{U_\delta}\eta w^2\,d\mu
    -\int_{U_\delta}|\nabla r|^2\,d\mu
    \geq \int_{U_\delta}\eta w^2\,d\mu-C\delta.
\end{align*}
By \eqref{eva_00} and $2\delta\leq\delta_0$, the source term satisfies
\begin{align*}
    \left|\int_B\eta\xi_\delta\frac{\rho}{w}
        \langle\nabla u,\nabla r\rangle\,d\mu\right|
    \leq \int_{U_{2\delta}}|\rho|\,dx.
\end{align*}
Finally, on $U_{2\delta}$ we have $r\geq1/2$, so the radial extension satisfies $|D\eta|\leq 2\norm{D^\top\eta}_{L^\infty(\del B)}$. Using $|\nabla\eta|\leq w|D\eta|$ and the third estimate in \eqref{eq:easy_energy_integral_bounds} with $2\delta$ in place of $\delta$, we obtain
\begin{align*}
    \left|\int_B\xi_\delta\langle\nabla r,\nabla\eta\rangle\,d\mu\right|
    \leq C\norm{D^\top\eta}_{L^\infty(\del B)}\sqrt{\delta}.
\end{align*}
Combining the last three estimates gives \eqref{eq:ikari}.

\medskip
Throughout Steps 3, 4 and 5, let $0<\eps<\eps_0$, with $\eps_0$ sufficiently small.

\subsubsection*{\textbf{Step 3. Smooth competitors via convolution}} \textit{Let $\varphi_\eps \in C_c^\infty(SO(m))$ be as in \eqref{eq:rotation_profile}.
% be a family of smooth, non-negative functions, such that \oorange{scegliere esplicitamente il profilo}
% \begin{align}\label{eq:properties_of_varphi_eps}
%     \operatorname{supp} \varphi_\eps \subseteq B_\eps^{SO(m)}(\id), \qquad \int_{SO(m)} \varphi_\eps(Q) \, dQ = 1,
% \end{align}
% where $dQ$ is the (normalised) Haar measure of $SO(m)$. 
Then the convolution}
\begin{align}\label{eq:convolution_on_the_sphere}
    u_\eps(x) = (P_\eps u)(x) \doteq \int_{SO(m)} u(Qx) \, \varphi_\eps(Q) \, dQ \qquad x \in \overline B
\end{align}
\textit{has the following properties.}
\begin{enumerate}[label = (\roman*)]
    \item\label{item:gradient_of_convolution} \textit{The gradient of the convolution $u_\eps$ and its tangential component are given by}
    \begin{align*}
        Du_\eps = \int_{SO(m)} Q^t Du(Q\cdot) \, \varphi_\eps(Q) \, dQ, \qquad 
        D^\top u_\eps = \int_{SO(m)} Q^t (D^\top u)(Q\cdot) \, \varphi_\eps(Q) \, dQ.
    \end{align*}
    \textit{In particular}
    \begin{align*}
        |Du_\eps|^2 \leq P_\eps\left(|Du|^2\right), \qquad |D^\top u_\eps|^2 \leq P_\eps\left( |D^\top u|^2 \right).
    \end{align*}
    \item\label{item:spacelike_convolutions} \textit{Each $u_\eps$ is strictly space-like and $\norm{u - u_\eps}_{L^\infty(B)}\to 0$ as $\eps\to 0$.}
    \item\label{item:area_convolutions} \textit{The convolution increases the area of the graph, that is $J(u_\eps) \leq J(u)$, where $J$ is defined in \eqref{eq:J}.}
    \item \textit{It holds}
    \begin{align}\label{eq:tangent_mollified_gradient}
            \norm{D^\top D^\top u_\eps}_{L^\infty(\del B)} \leq \frac{C}{\eps}, \qquad |D^\top u_\eps (r\omega) - D^\top u_\eps(\omega)| \leq C\frac{1-r}{\eps}.
    \end{align}
    \item\label{item:collar_convolutions} For every $\delta>0$ it holds
    \begin{align}\label{eq:collar_convolutions}
        \int_{U_\delta} w \left( 1 - Du \cdot Du_\eps \right) \, dx \leq C(\eps + \delta).
    \end{align}
\end{enumerate}

Since $D(u\circ Q) = Q^t \circ Du \circ Q$, the first identity in \ref{item:gradient_of_convolution} follows simply by differentiating under the integral sign. {The second follows by projecting the first in the tangential direction.} By Jensen's inequality, which applies, as $\varphi_\eps  dQ$ is a probability measure on $SO(m)$, we have
\begin{align*}
    |D u_\eps|^2 \leq \int_{SO(m)} |Du (Q\cdot)|^2 \, \varphi_\eps(Q) \, dQ = P_\eps\left( |Du|^2\right).
\end{align*}
The other inequality is analogous, and strict space-likeness of $u_\eps$ is immediate:
\begin{align*}
    |Du_\eps|^2 \leq P_\eps\left(|Du|^2 \right) \leq \norm{Du}_{L^\infty(B)}^2 < 1.
\end{align*}
To conclude the proof of \ref{item:spacelike_convolutions}, observe that by \eqref{eq:rotation_group_bounds} and the space-likeness of $u$ we have 
\begin{align*}
    |u_\eps(x) - u(x)| &\leq \int_{SO(m)} |u(Qx) - u(x)| \, \varphi_\eps(Q) \, dQ \leq \int_{SO(m)} |Qx - \id x| \, \varphi_\eps(Q) \, dQ \leq \eps.
\end{align*}
Item \ref{item:area_convolutions} follows again by Jensen's inequality and Fubini's Theorem.

For fixed, unit, orthogonal vectors $\omega, \tau$, let $A \in \mathfrak{so}(m)$ be such that
\begin{align*}
    A\omega = \tau, \qquad A\tau = - \omega, \qquad Ae = 0, \qquad e \perp \operatorname{span}\{\omega,\tau\}
\end{align*}
Since the derivative of the curve $t \mapsto e^{tA}r\omega$ at $t=0$ is $r\tau$, differentiating we obtain
\begin{align}\label{eq:eva03}
    D^\top u_\eps (r\omega) \cdot r\tau = \frac{d}{dt}\bigg\vert_{t=0} u_\eps(e^{tA}r\omega) = - \int_{SO(m)} u(Qr\omega) \, (X_A\varphi_\eps)(Q) \, dQ
\end{align}
where $X_A$ is the infinitesimal generator of $A$. Taking $r = 1$ and differentiating again we have
\begin{align*}
    \frac{d^2}{dt^2}\bigg\vert_{t=0} u_\eps(e^{tA}\omega) = - \int_{SO(m)} Q\tau \cdot Du(Q\omega) \, (X_A\varphi_\eps)(Q) \, dQ.
\end{align*}
Since the derivative of $\varphi_\eps$ is $O(\eps^{-n-1})$, where $n = \dim SO(m)$, and $|Q\tau \cdot Du(Q\omega)|\leq 1$,
\begin{align*}
    |D^\top D^\top u_\eps (\tau,\tau)| \leq \frac{C}{\eps}, \qquad \text{for every $\tau \in T_\omega \del B$ with $|\tau|=1$.}
\end{align*}
The first in \eqref{eq:tangent_mollified_gradient} follows by the arbitrariness of $\tau$.

For the second one, use twice \eqref{eq:eva03} to obtain
\begin{align*}
    \left(rD^\top u_\eps(r\omega) - D^\top u_\eps(\omega) \right) \cdot \tau &= \int_{SO(m)} \left(u(Q\omega)- u(Qr\omega) \right) (X_A\varphi_\eps)(Q) \, dQ \leq C\frac{1-r}{\eps}.
\end{align*}
Taking the supremum in $\tau \in T_\omega\del B$ with $|\tau|=1$ we obtain
\begin{align*}
    |rD^\top u_\eps(r\omega) - D^\top u_\eps(\omega)| \leq C \frac{1-r}{\eps}.
\end{align*}
To conclude observe that
\begin{align*}
    |D^\top u_\eps(r\omega) - D^\top u_\eps(\omega)| &\leq |D^\top u_\eps (r\omega) - rD^\top u_\eps(r\omega)| + |rD^\top u_\eps(r\omega) - D^\top u_\eps(\omega)| \\ &\leq (1-r)|D^\top u_\eps (\omega r)| + C\frac{1-r}{\eps} \leq C'\frac{1 - r}{\eps}.
\end{align*}

We finally prove \ref{item:collar_convolutions}. By minimality of $u$ and \ref{item:area_convolutions} and \ref{item:spacelike_convolutions}, we have
\begin{align*}
    0 \leq I(u_\eps) - I(u) 
    % = \int_B \left\{ \mathscr L(Du_\eps) - \mathscr L(Du) \right\} \, dx - \int_B (u_\eps-u) \, d\rho + \int_{\del B} (u_\eps-u), \, d\psi
    \leq \eps \left( \int_B d|\rho| + \int_{\del B} d|\psi| \right).
\end{align*}
Notice that
\begin{align*}
    I(u_\eps) - I(u) &= \int_B \left\{\mathscr L(Du_\eps) - \mathscr L (Du) - D\mathscr L(Du) \cdot(Du_\eps - Du) \right\} \, dx\\ &= \int_B\frac {w_u}{2} \left\{ |Du - Du_\eps|^2 + (\mathscr L(Du) - \mathscr L(Du_\eps))^2 \right\} \, dx \\ &= \int_B w(1 - Du\cdot Du_\eps) \, dx -  \int_B \sqrt{1-|Du_\eps|^2} \, dx.
\end{align*}
Since the integrand in the sum-of-squares expression is nonnegative, we can restrict the integral to $U_\delta$ for any $\delta\in (0,1)$. Since
\begin{align*}
    \int_{U_\delta} \sqrt{1 - |Du_\eps|^2} \, dx \leq |U_\delta| \leq C\delta,
\end{align*}
the thesis follows.

\subsubsection*{\textbf{Step 4. Angle estimates}\footnote{The first proof of this step contained a flaw pointed out and fixed by an AI tool, the conversation with which has also been fundamental in the writing of the {Appendix}.}} \textit{For a fixed $\theta \in (0,1)$ let}
\begin{align}
    E_\theta \doteq \set{\, |D^\top u_\eps| \geq 1 - \theta \, } \cap U_\delta, \qquad F_\theta \doteq \set{ \, |D^\top u_\eps|\geq 1 - \theta \, } \cap \del B.
\end{align}
\textit{There exists a constant $C = C(m,q)$ such that}
\begin{align}\label{eq:asuka}
    \int_{F_\theta} \cosh \gamma \, d\sigma \leq C \left( 1 + \norm{\psi}_{L^q(\del B)} \right) \theta^\beta, \qquad \beta \doteq \min \left\{\tfrac 12, 1 - \tfrac1q \right\},
\end{align}
\textit{uniformly in $0<\eps<\eps_0$.}

By \cref{lem:rotation_sphere_kernel}, the convolution can be also written as
\begin{align*}
    (P_\eps f )(\omega) = \int_{\del B} k_\eps(\omega,\omega') f(\omega') \, d\omega'
\end{align*}
for some kernel $\set{k_\eps}_{\eps}$ satisfying
\begin{align}\label{eq:katsuragi}
    \int_{\del B} k_\eps(\omega_0,\omega) \, d\omega = 1, \qquad \lambda \eps^{1-m} \bm 1_{\mathtt B_{\lambda\eps}(\omega_0)}(\omega) \leq k_\eps(\omega_0,\omega) \leq \Lambda \eps^{1-m} \bm 1_{\mathtt B_{\Lambda \eps}(\omega_0)}(\omega), \qquad \eps\ll1
\end{align}
where $\mathtt B_r (\omega_0)$ is the geodesic ball in $\del B$ of radius $r$ and centre $\omega_0 \in \del B$, for some uniform $0<\lambda\leq \Lambda$.

Recall that $\cosh\gamma \, d\sigma = w_\top^{-1}\cosh \gamma \, d\omega$. Observe that, for every $\omega_0 \in F_\theta$ it holds
\begin{align}\label{eq:rei}
    \int_{\mathtt B_{\lambda\eps}(\omega_0)} w_\top^{-2} \, d\omega \leq \frac{2\theta \eps^{m-1}}{\lambda}.
\end{align}
Indeed, by \ref{item:gradient_of_convolution}, we have
\begin{align}\label{eq:misato}
    P_\eps w_\top^{-2} = 1 - P_\eps \left(|D^\top u|^2 \right) \leq 2\theta \qquad \text{in }F_\theta.
\end{align}
% Hence, by the H\"older inequality
% \begin{align*}
%     P_\eps \left( w_\top^{-1} \, \cosh\gamma \right) &\leq \left( P_\eps \left( w_\top^{-q'} \right) \right)^{\frac{1}{q'}} \left( P_\eps \cosh^q\gamma \right)^{\frac{1}{q}} \\&\leq\left( P_\eps \left( w_\top^{-2} \right) \right)^{\frac{1}{q'}} \left( P_\eps \cosh^q\gamma \right)^{\frac{1}{q}}  \leq C\theta^{\frac{1}{q'}} \left( P_\eps \cosh^q\gamma \right)^{\frac{1}{q}} \qquad \text{in }F_\theta,
% \end{align*}
% where we used $q\leq 2$ and \eqref{eq:misato}.
Finally, by \eqref{eq:katsuragi}, for every $\omega_0 \in F_\theta$ we have
\begin{align*}
    \lambda \eps^{1-m} \int_{\mathtt B_{\lambda\eps}(\omega_0)} w_\top^{-2} \, d\omega \leq \left(P_\eps w_\top^{-2}\right)(\omega_0) \leq 2\theta,
\end{align*}
which is \eqref{eq:rei}.

Choose centres $\omega_1, \dots, \omega_h \in F_\theta$ such that
\begin{align*}
    F_\theta\subseteq \bigcup_{j=1}^h \mathtt B_{\lambda \eps} (\omega_j), \qquad h\eps^{m-1} \leq C, \qquad \sum_{j=1}^h \bm 1_{\mathtt B_{\lambda\eps}(\omega_j)} \leq C.
\end{align*}
Applying \eqref{eq:rei} to each ball and summing up we obtain
\begin{align*}
    \int_{F_\theta} w_\top^{-2} \, d\omega \leq C\theta.
\end{align*}
Now let $s \doteq \min\set{2,q}$, so that $s'\geq 2$. By H\"older's inequality
\begin{align*}
    \int_{F_\theta} \cosh \gamma \, d\sigma \leq \norm{\cosh \gamma}_{L^s(\del B)} \left(\int_{F_\theta} w_\top^{-s'} \, d\omega \right)^{\frac{1}{s'}} \leq C (1 + \norm{\psi}_{L^q(\del B)}) \theta^{\frac{1}{s'}}
\end{align*}

\subsubsection*{\textbf{Step 5. Conclusion}}
Assume that $\theta<1/3$ and $C\delta/\eps\leq\theta$, where $C$ is the constant in \eqref{eq:tangent_mollified_gradient}. Choose a cut-off $\eta\in C^\infty(\del B)$ such that
\begin{align*}
    0\leq\eta\leq 1, \qquad
    \eta=1 \quad \text{on }F_{2\theta}, \qquad
    \eta=0 \quad \text{on }F_{3\theta}^c,
\end{align*}
and extend it radially to $U_{2\delta}$. Thanks to the first estimate in \eqref{eq:tangent_mollified_gradient}, we can also require
\begin{align*}
    \norm{D^\top\eta}_{L^\infty(\del B)}
    \leq \frac{C}{\theta\eps}.
\end{align*}
Moreover, if $r\omega\in E_\theta$, the second estimate in \eqref{eq:tangent_mollified_gradient} gives
\begin{align*}
    |D^\top u_\eps(\omega)|
    \geq |D^\top u_\eps(r\omega)|-C\frac{1-r}{\eps}
    \geq 1-\theta-C\frac{\delta}{\eps}
    \geq 1-2\theta.
\end{align*}
Thus $\omega\in F_{2\theta}$ and $\eta=1$ on $E_\theta$. By combining \eqref{eq:ikari} and \eqref{eq:asuka}, we obtain
\begin{align}
    \int_{E_\theta}w^2\,d\mu
    &\leq \int_{U_\delta}\eta w^2\,d\mu \nonumber\\
    &\leq C\delta
        +C\norm{D^\top\eta}_{L^\infty(\del B)}
        \sqrt{\delta}
        +\int_{\del B}\eta\cosh\gamma\,d\sigma + \int_{U_{2\delta}} |\rho| \, dx \nonumber\\
    &\leq C\delta
        +\frac{C}{\theta\eps}\sqrt{\delta}
        +\int_{F_{3\theta}}\cosh\gamma\,d\sigma + \int_{U_{2\delta}} |\rho| \, dx \nonumber\\
    &\leq C\left(\delta+
        \frac{\sqrt{\delta}}{\theta\eps}+\theta^\beta\right) + \int_{U_{2\delta}} |\rho| \, dx.
        \label{eq:eva02}
\end{align}

On the other hand, \eqref{eq:easy_energy_integral_bounds} and \eqref{eq:collar_convolutions} control the energy in $E_\theta^c=U_\delta\setminus E_\theta$:
\begin{align}\label{eq:fujutsuki}
    \int_{E_\theta^c}w\,dx
    \leq \frac{C}{\theta}\left(\eps+\sqrt{\delta}\right).
\end{align}
Indeed, on $E_\theta^c$ we have
\begin{align*}
    Du\cdot Du_\eps
    &=D^\top u\cdot D^\top u_\eps+\del_r u\,\del_r u_\eps\\
    &\leq |D^\top u_\eps|+|\del_r u|
    \leq 1-\theta+|\del_r u|.
\end{align*}
Hence, by \eqref{eq:collar_convolutions},
\begin{align*}
    \theta\int_{E_\theta^c}w\,dx
    &\leq \int_{E_\theta^c}w(1-Du\cdot Du_\eps)\,dx
        +\int_{E_\theta^c}w|\del_r u|\,dx\\
    &\leq C(\eps+\delta)+\int_{E_\theta^c}w|\del_r u|\,dx.
\end{align*}
Finally, since $\langle\nabla u,\nabla r\rangle=w^2\del_r u$ and $d\mu=w^{-1}\,dx$, by \eqref{eq:easy_energy_integral_bounds} we have
\begin{align*}
    \int_{E_\theta^c}w|\del_r u|\,dx
    &\leq \int_{U_\delta}|\langle\nabla u,\nabla r\rangle|\,d\mu\\
    &\leq \left(\int_{U_\delta}w^2\,d\mu\right)^{1/2}
        \left(\int_{U_\delta}|\nabla r|^2\,d\mu\right)^{1/2}
    \leq C\sqrt{\delta},
\end{align*}
and \eqref{eq:fujutsuki} follows.

Summing \eqref{eq:eva02} and \eqref{eq:fujutsuki}, and using $w^2\,d\mu=w\,dx$, we deduce
\begin{align*}
    \int_{U_\delta}w\,dx
    \leq C\left(\delta+
        \frac{\sqrt{\delta}}{\theta\eps}
        +\theta^\beta
        +\frac{\eps+\sqrt{\delta}}{\theta}\right)+
        \int_{U_{2\delta}} |\rho| \, dx.
\end{align*}
Choose
\begin{align*}
    \eps=\delta^{1/4}, \qquad
    \theta=\delta^{1/[4(\beta+1)]}.
\end{align*}
For $\delta$ sufficiently small, this choice satisfies $\theta<1/3$ and $C\delta/\eps\leq\theta$. Since
\begin{align*}
    \frac{\sqrt{\delta}}{\theta\eps}
    =\frac{\eps}{\theta}
    =\theta^\beta
    =\delta^{\beta/[4(\beta+1)]},
\end{align*}
we obtain the desired estimate with $\alpha_q=\beta/[4(\beta+1)]$. For larger $\delta$, the estimate follows from the global bound in \eqref{eq:easy_energy_integral_bounds} after increasing $C$. This concludes the proof of \cref{thm:boundary_energy_estimate}. \qed

\section{Proof of the main results}

We construct the solution by approximation with solutions corresponding to smoothed data. Hence, a preliminary step is the following existence result.

\begin{prop}\label{prop:smooth_capillary_existence}
    Let $\rho \in C^\infty(\overline B)$ and $\psi \in C^\infty(\del B)$. Then the problem \eqref{eq:problem_Rm} has a unique zero-mean, smooth, strictly space-like, classical solution $u \in C^\infty(\overline B)$ if and only if $\rho$ and $\psi$ meet the compatibility condition \eqref{eq:compatibility}. Moreover
    \begin{align*}
        \sup w_u \leq C \left(m,R, \norm{\rho}_{L^\infty(B)}, \norm{\psi}_{\lip(\del B)} \right).
    \end{align*}
\end{prop}
\begin{proof}
    Necessity of the compatibility condition is clear. Let us assume it and prove existence and regularity of a solution. We will do it by showing that the set $\mathscr T$ of $t \in [0,1]$ such that the problem
    \begin{align}\label{eq:t_problem}
        -\div\Bigg(\dfrac{Du}{\sqrt{1 - |Du|^2}}\Bigg)= t\rho \quad \text{in $B$,} \qquad -\dfrac{Du \cdot n}{\sqrt{1 - |Du|^2}} = t\psi \quad \text{on $\del B$}, \qquad \int_B u \, dx = 0,
    \end{align}
    has a solution $u_t$ with the desired properties is open and closed in $[0,1]$. Since $0 \in \mathscr T$, by connectedness this implies $\mathscr T = [0,1]$, which yields the conclusion.

    Openness follows by the implicit function Theorem. Fix $\alpha \in (0,1)$ and consider the Banach spaces
    \begin{align*}
        X \doteq  \left\{ \, v\in C^{2,\alpha}(\overline B) \ | \ \int_B v \, dx = 0 \,  \right\}, \qquad Y \doteq \left\{ (f,g) \in C^{0,\alpha}(\overline B) \times C^{1,\alpha}(\del B) \ | \ \int_B f \, dx = \int_{\del B} g \, d\omega \right\}
    \end{align*}
    and the map
    \begin{align*}
        \mathscr F: \mathscr U \times [0,1] \longrightarrow Y, \qquad \mathscr F(u,t) = (- \div(w_u Du) - t\rho, -w_u Du \cdot n - t\psi),
    \end{align*}
    where $\mathscr U \subseteq X$ is the open set of strictly space-like functions. If $t_0 \in \mathscr T$, and $u_0$ solves \eqref{eq:t_problem} with $t = t_0$, then $\mathscr F$ vanishes in $(u_0,t_0)$ and its differential with respect to the first variable is an isomorphism by standard elliptic theory (the kernel is the space of constants, which is trivialised by the zero-mean condition). The implicit function Theorem gives a continuous function from a neighbourhood of $t_0$ into $\mathscr U$ that maps $t$ to the solution $u_t$ of \eqref{eq:t_problem}. This shows that $\mathscr T$ is open.

    We are left with proving that $\mathscr T$ is closed. Notice first that every solution entering the definition of $\mathscr T$ is covered by \cref{thm:full_gradient_estimate}, with the same bounds
    \[
    \Lambda=\norm{\rho}_{L^\infty(B)},\qquad \Psi_0=\norm{\psi}_{L^\infty(\partial B)},\qquad \Psi_1=\norm{D^\top\psi}_{L^\infty(\partial B)}.
    \]
    Hence, for every $t \in \mathscr T$ the solutions $u_t$ to \eqref{eq:t_problem} all enjoy the gradient bound $\norm{w_t}_{L^\infty(B)}\leq C$ uniformly in $t\in \mathscr T$, where $w_t = (1 - |Du_t|^2)^{-1/2}$. Hence, they are solutions to elliptic problems with uniformly controlled ellipticity constants and Lieberman's theory \cite[Theorem 2]{Lieberman88} applies yielding a uniform $C^{1,\beta}(\overline B)$ bound
    \begin{align*}
        \norm{u_t}_{C^{1,\beta}(\overline B)} \leq K, \qquad \text{for every $t \in \mathscr T$.}
    \end{align*}

    Take now a sequence $t_j \in \mathscr T$ converging to $t \in [0,1]$. Uniform spacelikeness and the H\"older estimate imply that, up to a subsequence, the solutions $u_j = u_{t_j}$ converge in $C^{1,\beta'}(\overline B)$, $0<\beta'<\beta$, to a space-like function $u \in C^{1,\beta'}(\overline B)$ with $\norm{w_u}_{L^\infty(B)} \leq C$. Hence we can pass to the limit in the weak solution identity
    \begin{align*}
        \int_B w_j Du_j \cdot D\eta \, dx = t_j \int_B \eta \rho \, dx - t_j \int_{\del B} \eta \psi \, d\omega, \qquad \text{for every $\eta \in C^1(\overline{B})$.}
    \end{align*}
    Thus $u$ weakly solves \eqref{eq:t_problem} with $t=\lim_jt_j$. Standard regularity theory gives smoothness up to the boundary, so the equation holds pointwise.
\end{proof}

\subsection{Reduction to uniform integrability of $\set{w_j}$}\label{subsec:mollify} Now let $\rho \in \mathscr M(B)$ and $\psi \in L^q(\del B)$, $1<q\leq \infty$, satisfy the compatibility condition \eqref{eq:compatibility}. Let $\tilde \rho \in \mathscr M(\R^m)$ be the extension of $\rho$ that vanishes outside $B$ and consider the mollification
\begin{align*}
    \tilde \rho_j \doteq \tilde \rho * \zeta_{1/j}
\end{align*}
where $\{\zeta_\eps\}_{\eps>0}$ is a standard regularising kernel in $\R^m$. We also regularise $\psi$ by convolution on $\del B$
\begin{align*}
    \psi_j \doteq P_{1/j}\psi,
\end{align*}
where $P_\eps$ is defined in \eqref{eq:convolution_on_the_sphere}. In order for the compatibility condition to be met at every $j$ we modify $\tilde \rho_j$ by
\begin{align*}
    \rho_j \doteq \tilde\rho_j - c_j, \qquad
    c_j \doteq \frac{1}{|B|} \left( \int_B \tilde\rho_j \, dx - \int_{\del B} \psi_j \, d\omega \right) .
\end{align*}
By standard theory
\begin{gather}\label{eq:general_smoothened}
    \rho_j \in C^\infty(\overline B), \qquad \psi_j \in C^\infty(\del B),\qquad c_j \to 0, \nonumber \\ \qquad \rho_j \, dx \overset{*}{\rightharpoonup} \rho \quad \text{in $\mathscr M(\overline B)$}, \qquad \psi_j \to \psi \quad \text{strongly in $L^q(\del B)$ for $q<\infty$ and weakly-* in $L^\infty(\del B)$.}
\end{gather}
Let $u_j \in \mathscr A(B)$ be the zero-mean solution given by \cref{prop:smooth_capillary_existence}, and write $w_j = (1 - |Du_j|^2)^{-1/2}$. By \eqref{eq:balthasar} in \cref{lem:magi}, $u_j$ is the minimiser of $I_j = I_{\rho_j,\psi_j}$ in $\mathscr A(B)$. By \eqref{eq:casper} in \cref{lem:magi}, up to taking a subsequence, we have
\begin{align*}
    u_j \to u, \qquad \text{uniformly in $\overline B$}, \qquad Du_j \to Du \quad \text{a.e. in $B$.}
\end{align*}
In order to conclude that $u$ is a weak solution to \eqref{eq:problem_Rm} in a relatively open set $\Omega \subseteq \overline{B}$, we must pass to the limit in both sides of
\begin{align}\label{eq:approximisers}
    \int_B w_j Du_j \cdot D\eta \, dx = \int_B \eta \rho_j \, dx - \int_{\del B} \eta\psi_j \, d\omega \qquad \text{for every $\eta \in \lip_c(\Omega)$}.
\end{align}
By \eqref{eq:general_smoothened}, it suffices to show that
\begin{align}\label{eq:the_goal}
    \lim_{j\to \infty} \int_B w_j Du_j \cdot D\eta \, dx = \int_B w Du \cdot D\eta \, dx \qquad \text{for every $\eta \in \lip_c(\Omega)$}.
\end{align}
By Vitali's Theorem, \eqref{eq:the_goal} holds true provided that $\set{w_j}$ is locally uniformly integrable in $\Omega$.
% \begin{align}\label{eq:vitali}
%     \text{\eqref{eq:the_goal}} \qquad \Longleftarrow \qquad \text{$\set{w_j}$ is locally uniformly integrable in $\Omega$.}
% \end{align}

\subsection{Proof of \cref{thm:existence_lip}}\label{subsec:proof_1}
If $\rho \in L^\infty(B)$ and $\psi \in \lip(\del B)$, by convolution we construct approximating sequences
\begin{align*}
    \rho_j \rightharpoonup \rho \quad \text{in $L^2(B)$}, \qquad \psi_j \to \psi\quad  \text{in }C(\del B) \cap W^{1,q}(\del B) \text{ for every $1\leq q <\infty$.}
\end{align*}
Moreover
\begin{align}\label{eq:norms_convolutions_bounded}
    \norm{\rho_j}_{L^\infty(B)} \leq 3 \norm{\rho}_{L^\infty(B)}, \qquad \norm{\psi_j}_{L^\infty(\del B)} \leq \norm{\psi}_{L^\infty(\del B)}, \qquad \norm{D^\top \psi_j}_{L^\infty(\del B)} \leq \norm{D^\top \psi}_{L^\infty(\del B)}.
\end{align}
For every $j$, let $u_j$ be the smooth, zero-mean solution to \eqref{eq:problem_Rm} with data $\rho_j$ and $\psi_j$ established by \cref{prop:smooth_capillary_existence}. By \cref{thm:full_gradient_estimate} a gradient bound $\norm{w_j}_{L^\infty(B)} \leq C$ holds along the sequence, uniformly in $j$. Lieberman's result \cite[Theorem 2]{Lieberman88} hence yields a H\"older bound $\norm{u_j}_{C^{1,\beta}(\overline B)} \leq K$ uniformly in $j$. Hence, up to a subsequence
\begin{align}\label{eq:djo}
    u_j \to u \quad \text{in $C^{1,\beta'}(\overline B)$,}
\end{align}
for some $\beta'<\beta$. By Lebesgue's dominated convergence, \eqref{eq:the_goal} holds with $\Omega = \overline B$ and this makes of $u$ a weak solution to \eqref{eq:problem_Rm} on the whole of $\overline B$. On the other hand, \cref{lem:L^2_w_estimates} implies that $\set{u_j}$ is bounded in $W^{2,2}(B)$.
% and more precisely
% \begin{equation}\label{eq:harry_styles}
%     \int_B |D^2u_j|^2 \, dx \leq \int_B \rho_j^2 \, dx + R \int_{\del B} |D^\top \psi_j|^2 \, d\omega.
% \end{equation}
Hence, up to a subsequence $u_j \rightharpoonup v$ in $W^{2,2}(B)$ and strongly in $W^{1,2}(B)$. Uniqueness of the limit and \eqref{eq:djo} imply that $u = v \in W^{2,2}(B)$. The integral estimate \eqref{eq:uniform_second_control} follows by \cref{lem:L^2_w_estimates} and \eqref{eq:dinner_party}.

% If now $u_1,u_2$ are two space-like, weak solutions, then testing the weak formulations with $u_2 - u_1$ and subtracting gives
% \begin{align*}
%     \int_B (w_{u_2} Du_2 - w_{u_1} Du_1 ) \cdot D(u_2 - u_1) \, dx = 0.
% \end{align*}
% By convexity, the integrand bounds $|Du_1 - Du_2|^2$, hence it must be $u_1 = u_2 + \mathrm{const.}$ a.e. in $B$. This concludes the proof.

\subsection{Proof of \cref{thm:existence_H^1}} If $\rho \in L^2(B)$ and $\psi \in W^{1,2}(\del B)$ then
\begin{align*}
    \rho_j \rightharpoonup \rho \quad \text{in $L^2(B)$}, \qquad \psi_j \to \psi\quad  \text{in } W^{1,2}(\del B).
\end{align*}
Moreover
\begin{align}\label{eq:dinner_party}
    \norm{\rho_j}_{L^2(B)} \leq 3 \norm{\rho}_{L^2(B)}, \qquad \norm{D^\top \psi_j}_{L^2(\del B)} \leq \norm{D^\top \psi}_{L^2(\del B)}.
\end{align}
As in the proof of \cref{thm:existence_lip}, \cref{subsec:proof_1}, we construct a sequence of smooth, space-like functions
\begin{align*}
    u_j \to u \quad \text{in $C(\overline B)$}
\end{align*}
that solve \eqref{eq:problem_Rm} with respect to the smoothened sources $\rho_j$, $\psi_j$. This time, however, uniform space-likeness is generally lost in the limit. We circumvent this thanks to our second fundamental form estimates. By \cref{lem:L^2_w_estimates} and \eqref{eq:dinner_party} we have
\begin{align*}
    \int_B |D^2u_j|^2 \, dx \leq 9 \int_B \rho^2 \, dx + R \int_{\del B} |D^\top \psi|^2 \, d\omega,
\end{align*}
hence $\set{u_j}$ is bounded in $W^{2,2}(B)$. Up to a subsequence
\begin{align*}
    u_j \rightharpoonup v \quad \text{in $W^{2,2}(B)$,} \qquad u_j \to v \quad \text{ in $W^{1,2}(B)$}, \qquad u_j \to v \quad \text{a.e. in $B$,} \qquad Du_j \to Dv \quad \text{a.e. in $B$.}
\end{align*}
By uniqueness of the limit, it must be $u = v \in W^{2,2}(B)$. Moreover, by \cref{lem:L^2_w_estimates} one has
\begin{align*}
    \sup _j\int_B w_j^2 \, dx < \infty
\end{align*}
hence, by de la Vallée-Poussin's Theorem, $\set{w_j}$ is locally uniformly integrable and \eqref{eq:the_goal} follows: $u$ weakly solves \eqref{eq:problem_Rm}.
% \begin{align*}
%     \lim_{j\to \infty} \int_B w_j Du_j \cdot D\eta \, dx = \int_B wDu \cdot D\eta \, dx, \qquad \text{for every $\eta \in C^1(\overline B)$.}
% \end{align*}
% Therefore, we can pass to the limit in the weak formulation and deduce that the limit $u$ is a weak solution to \eqref{eq:problem_Rm} with sources $\rho,\psi$.
The estimates \eqref{eq:uniform_second_control} and \eqref{eq:uniform_w_control} follow from \eqref{eq:dinner_party} and \cref{lem:L^2_w_estimates}.
\qed

\subsection{Proof of \texorpdfstring{\cref{thm:ultimate_thm}}{Theorem \ref{thm:ultimate_thm}}}

After rescaling we can assume $B = B_1$. Let $\rho = f \, dx + \nu$ and $\psi$ be as in the hypotheses of \cref{thm:ultimate_thm}. Extending $f$ to zero outside the ball, the approximation in \cref{subsec:mollify} satisfies
\begin{align*}
    \rho_j = f_j + \nu_j - c_j, \qquad
    f_j = f * \zeta_{1/j}, \qquad
    \nu_j = \nu * \zeta_{1/j},
\end{align*}
where as above $c_j$ are constants that ensure that the compatibility condition is met by $\rho_j$ and $\psi_j$ for every $j$. Notice that $c_j \to 0$ and
\begin{align}\label{eq:ultimate_smoothened_convergence}
    \rho_j \, dx \overset{*}{\rightharpoonup} \rho \quad \text{in $\mathscr M(B)$}, \qquad f_j \to f \quad \text{in $L^1(B)$}, \qquad \rho_j \to f \quad \text{in $L^2_{\loc}(B\setminus\Sigma)$}.
\end{align}
Let
\begin{align*}
    M \doteq \sup_j \left( \norm{\rho_j}_{L^1(B)} + \norm{\psi_j}_{L^1(\del B)} \right).
\end{align*}

We test the equation of $u_j$ with $u_j - u_j(0)$ to obtain
\begin{align*}
    \int_B w_j \, dx = \int_B w_j^{-1} \, dx + \int_B (u_j - u_j(0)) \rho_j \, dx - \int_{\del B} (u_j - u_j(0)) \psi_j \, d\omega \leq |B| + M
\end{align*}
for every $j$. Fatou's Lemma implies $w \in L^1(B)$ and, in particular, $|Du|<1$ a.e. in $B$. By \cite[Theorem 1.11]{DeFilippisMingione26}, this implies that, if $m =2$, then there are no interior light segments for the minimiser $u$, and, since $B$ is strictly convex, there are no light segments at all. This settles \ref{item:m=2}.

Write $K = K_\psi^\rho$ and let $\Omega''\Subset \Omega' \subseteq \Omega\Subset B\setminus (\Sigma\cup K)$ be nested subdomains. By \cite[Lemma 2.4]{Maniscalco26}, the restriction $u_j|_\Omega$ (resp. $u|_{\Omega}$) is the unique minimiser of the functional $I_{\rho_j,0}$ (resp. $I_{\rho,0}$) in the class of $1$-Lipschitz functions that equal $u_j$ (resp. $u$) on $\del\Omega$. Since $u$ does not have light segments in $B\setminus K_\psi^\rho$, and that every $u_j$ is strictly space-like, the class $\mathscr F \doteq \set{u_j,u}$ is a family of admissible functions (namely $\mathscr F \subseteq \mathcal S(\del \Omega)$, cf. \cite{BIMM}) that is compact in $C(\del \Omega)$. For every weakly space-like $v$, let
\begin{align*}
    L_R^v(x) \doteq \set{ \, y \in B \ | \ |x-y|^2 - (v(x) - v(y))^2 \leq R^2 \, }, \qquad L_R^v(A) \doteq \bigcup_{x \in A} L_R^v(x)
\end{align*}
be the \emph{Lorentzian ball} centred in $x$ or in a subset $A\subseteq B$. The fact that $u$ does not have a light segment implies that there exists an $R>0$ such that
\begin{align*}
    L_R^{u_j}(\Omega'') \Subset \Omega', \qquad L_R^u(\Omega'') \Subset \Omega'.
\end{align*}
Moreover \eqref{eq:ultimate_smoothened_convergence} implies
\begin{align*}
    \sup_j \int_{\Omega'} \rho_j^2 \, dx < \infty.
\end{align*}
By \cite[Theorem 4.13]{BIMM} we obtain
\begin{align*}
    \sup_j \int_{\Omega''} \left\{ w_j \ln w_j + |D^2u_j|^2 \right\} \,dx < \infty.
\end{align*}
By the arbitrariness of $\Omega''$ and $\Omega'$, this yields $u \in W_\loc^{2,2}(B\setminus(\Sigma \cup K))$. Hence, by de la Vallée-Poussin's Theorem, $\set{w_j}$ is locally uniformly integrable in $B\setminus (\Sigma\cup K)$. By \cite[Theorem 4.2]{BIMM} and $\Haus^1(\Sigma) = 0$, we conclude that $\set{w_j}$ is locally uniformly integrable in the whole of $B\setminus K$. We are left with showing uniform integrability up to the boundary.% i.e. that for every $\eps>0$ there exists $\alpha>0$ such that, for every measurable set $A\subseteq B\setminus K$ it holds
% \begin{align*}
%     |A| < \alpha \qquad \Longrightarrow \qquad \sup_j\int_A w_j \, dx < \eps.
% \end{align*}

By \cref{thm:boundary_energy_estimate} we have
\begin{align*}
    \int_{U_\delta} w_j \, dx \leq C \left( \delta^{\alpha_q} + \int_{U_{2\delta}} \left\{|f_j| + |c_j| \right\} \, dx \right)
\end{align*}
for every $j\gg1$ and $\delta>0$ small enough so that $\nu_j = 0$ on $U_{2\delta}$. Since $f_j \to f$ in $L^1(B)$ and $c_j \to 0$, this yields
\begin{align*}
    \lim_{\delta \to 0}\sup_j \int_{U_\delta} w_j \, dx = 0.
\end{align*}
Fix $\Omega'\subseteq \overline B\setminus K$ and $\eps>0$. We first fix $\delta>0$ small enough so that
\begin{align*}
    \sup_j\int_{U_\delta} w_j \, dx < \frac\eps2.
\end{align*}
Since $\overline \Omega'$ si a compact subseet of $B\setminus K$, local uniform integrability gives $\alpha>0$ such that, for every measurable $A\subseteq \Omega'$ with $|A|<\alpha$, %By local uniform integrability in $B\setminus K$, for every $\Omega\Subset B_{1-\delta}\setminus K$ there exists $\alpha = \alpha(\eps,\Omega)>0$ such that, for any measurable $A\subseteq \Omega$, if $|A\cap \overline{B_{1-\delta}}|<\alpha$ then
\begin{align*}
    \sup_j \int_A w_j \, dx = \sup_j \int_{U_\delta} w_j \, dx + \sup_j \int_{A\cap \overline{B_{1-\delta}}} w_j \, dx < \frac{\eps}{2} + \frac\eps2 = \eps.
\end{align*}
That is, $\set{w_j}$ is locally uniformly integrable on the whole of $\overline B\setminus K$. By \cref{eq:the_goal}, $u$ satisfies the weak formulation there.

% As a consequence, by Vitali's convergence Theorem
% \begin{align*}
%     w_j \to w, \qquad w_j Du_j \to wDu \qquad \text{in $L^1_{\loc}(B\setminus \Sigma)$.}
% \end{align*}
% We can thus pass to the limit in the weak formulation for each $u_j$, obtaining
% \begin{align*}
%     \int_B w Du \cdot D\eta \, dx = \int_B \eta \, d\rho - \int_{\del B} \eta \, d\psi, \qquad \text{for every $\eta \in \lip_c(B\setminus K)$,}
% \end{align*}
% as desired.

Finally, let $\Omega'\subseteq B\setminus(\Sigma \cup K)$ be open with $\nu|_{\Omega'}=0$ and $f\in L^\infty_{\loc}(U)$. On each $\Omega'' \Subset \Omega'$, $u$ minimises the Dirichlet functional with source $f$ and its own trace, by \cite[Lemma 2.4]{Maniscalco26}. Its trace on $\del \Omega'$ is space-like. The interior regularity result in \cite{BIMM} gives $u\in C^{1,\alpha}_{\loc}(\Omega')$ and $|Du|<1$ in $\Omega'$. If $f$ is smooth in $\Omega'$, local elliptic regularity gives $u\in C^\infty(\Omega')$. 

Uniqueness up to additive constant follows from \eqref{eq:melchior} and \eqref{eq:balthasar} of \cref{lem:magi} after subtracting the mean.\qed

\appendix

\section{Regularisation by rotations}\label{app:rotational_convolution}
Let $n=m(m-1)/2$, the dimension of $SO(m)$. We give $SO(m)$ the bi-invariant metric induced by $\langle A,B\rangle=-\tfrac12\tr(AB)$ on $\mathfrak{so}(m)$ and denote its distance by $\de$. The Haar measure $dQ$ on $SO(m)$ has mass one. On $\del B_1$ we use the usual surface measure $d\omega$.

\emph{The profile on the group.}
Fix $\eps_0<1$ below the injectivity radius of $SO(m)$. For $0<\eps<\eps_0$ set
\begin{align*}
  b(t)=
  \begin{cases}
    \exp\left(-\tfrac{1}{1-4t}\right) & 0\leq t<1/4,\\
    0 & t\geq 1/4 ,
  \end{cases}
  \qquad
  Z_\eps=\int_{SO(m)} b\left(\frac{\de(\id,Q)^2}{\eps^2}\right)\,dQ .
\end{align*}
Define the regularising kernel $\set{\varphi_\eps}_{\eps\in(0,\eps_0)}$ by
\begin{equation}\label{eq:rotation_profile}
   \varphi_\eps(Q)=\frac{1}{Z_\eps} b\left(\frac{\de(\id,Q)^2}{\eps^2}\right).
\end{equation}
Each $\varphi_\eps$ is smooth, non-negative, supported in $B_\eps^{SO(m)}(\id)$ and has integral one. It is invariant under inversion and under conjugation by every $T\in O(m)$. In exponential coordinates $Q=\exp A$ near $\id$, we have $\de(\id,Q)=|A|$ and Haar measure has a smooth positive density. Consequently $c\eps^n\leq Z_\eps\leq C\eps^n$ and
\begin{equation}\label{eq:rotation_group_bounds}
   c\eps^{-n}\bm 1_{B^{SO(m)}_{\eps/4}(\id)}
   \leq\varphi_\eps\leq C\eps^{-n}\bm 1_{B^{SO(m)}_{\eps/2}(\id)} .
\end{equation}
For $A\in\mathfrak{so}(m)$ let $X_A$ be its infinitesimal generator, i.e. $X_A h(Q)=\left.\tfrac{d}{dt}\right|_{t=0}h(Qe^{tA})$ for every smooth function $h$. In the same coordinates, differentiating the rescaled bump gives, for $|A_i|\leq1$ and $\ell\geq1$,
\begin{equation}\label{eq:rotation_profile_derivatives}
   \|X_{A_1}\cdots X_{A_\ell}\varphi_\eps\|_\infty
       \leq C_\ell\eps^{-n-\ell},\qquad
   \|X_{A_1}\cdots X_{A_\ell}\varphi_\eps\|_1
       \leq C_\ell\eps^{-\ell}.
\end{equation}
Indeed, the coefficients of the invariant fields and their derivatives are bounded in a fixed coordinate neighbourhood of the identity, and $\spt\varphi_\eps$ has Haar measure $O(\eps^n)$. For the general construction of a smooth approximate identity in exponential coordinates see the discussion following \cite[Definition~2.7]{AntonelliLeDonne2022}.

\begin{lem}\label{lem:rotation_sphere_kernel}
The operator $P_\eps f(\omega)=\int_{SO(m)} f(Q\omega)\varphi_\eps(Q)\,dQ$ has a smooth non-negative kernel with respect to $d\omega$,
\begin{align*}
    P_\eps f(\omega)=\int_{\del B} k_\eps(\omega,\omega')f(\omega')\,d\omega',
    \qquad \int_{\del B} k_\eps(\omega,\omega')\,d\omega'=1 .
\end{align*}
The kernel depends only on $\omega\cdot\omega'$. There are constants $0<\lambda\leq\Lambda$, depending only on $m$, such that for all sufficiently small $\eps$,
\begin{equation}\label{eq:rotation_sphere_bounds}
  \lambda\eps^{1-m}\bm 1_{\mathtt B_{\lambda\eps}(\omega)}(\omega')
  \leq k_\eps(\omega,\omega')
  \leq \Lambda\eps^{1-m}\bm 1_{\mathtt B_{\Lambda\eps}(\omega)}(\omega') ,
\end{equation}
where $\mathtt B_s(\omega)$ is a geodesic ball in $\del B$.
\end{lem}

\begin{proof}
Fix $e\in \del B$, put $H=\operatorname{Stab}(e)$, and write $\omega=ge$, $\omega'=g'e$. By Weil's integration formula
\cite[Theorem~2.51, formula~(2.52)]{Folland2015},
\begin{align*}
  k_\eps(\omega,\omega')
  =\frac{1}{|\del B|}\int_H\varphi_\eps(g'hg^{-1})\,dh ,
\end{align*}
where $dh$ is the Haar probability measure on $H$. Smooth local choices of $g,g'$ give smoothness, and $P_\eps 1=1$ gives the normalisation. The $O(m)$-conjugation invariance of $\varphi_\eps$ makes $k_\eps$ invariant under the diagonal action of $O(m)$, so it depends only on $\omega\cdot\omega'$.

Let $s=\dist_{\del B}(\omega,\omega')$. Since $Q\mapsto Q\omega$ is $1$-Lipschitz, the kernel vanishes for $s>\eps/2$. Otherwise choose a rotation $Q_0$ taking $\omega$ to $\omega'$ with $\de(\id,Q_0)=s$. Writing the fibre as $Q_0H_\omega$, the triangle inequality gives
\begin{align*}
  Q_0h\in B^{SO(m)}_{\eps/2}(\id)\ \Longrightarrow\ h\in B^{SO(m)}_\eps(\id),
  \qquad
  s<\eps/8,\ h\in B^{SO(m)}_{\eps/8}(\id)\
  \Longrightarrow\ Q_0h\in B^{SO(m)}_{\eps/4}(\id).
\end{align*}
Since $H_\omega$ is conjugate to $SO(m-1)$, the Haar probability of $H_\omega\cap B^{SO(m)}_t(\id)$ is comparable to $t^{n-(m-1)}$ for small $t$, uniformly in $\omega$. The two bounds now follow from \eqref{eq:rotation_group_bounds}.
\end{proof}

\bigskip

\noindent\hfill
\begin{tabular}{@{}l@{}}
Dipartimento di Matematica ``Giuseppe Peano''\\
Universit\`a degli Studi di Torino\\
Via Carlo Alberto 10, 10123 Torino, Italy\\
\texttt{lorenzo.maniscalco@unito.it}
\end{tabular}

% \addcontentsline{toc}{section}{Bibliography}
% \printbibliography[title={Bibliography}]

\end{document}